\documentclass[reqno,a4]{amsart}
\usepackage{amsmath, amssymb, amsthm}
\usepackage[dvipdfmx]{graphicx}

\usepackage{here}
\usepackage{url}

\usepackage[foot]{amsaddr}

\usepackage[monochrome]{color}

\allowdisplaybreaks

\newtheorem{theorem}{Theorem}
\newtheorem{defn}{Definition}
\newtheorem{proposition}{Proposition}
\newtheorem{lemma}{Lemma}
\newtheorem{remark}{Remark}

\begin{document}

\title{Non-self-intersective dragon curves}
\author{Shigeki Akiyama}
\address{Institute of Mathematics, University of Tsukuba, Japan}
\email{akiyama@math.tsukuba.ac.jp}
\author{Yuichi Kamiya}
\address{Department of Modern Economics, Faculty of Economics, Daito Bunka University, Japan}
\email{ykamiya@ic.daito.ac.jp}
\author{Fan Wen}
\address{Department of Mathematics, Jinan University, Guangzhou, China}
\email{wenfan\_fan@163.com (Corresponding author)}
\date{}
\maketitle

\begin{abstract}
Let us fold a strip of paper many times in the same direction,
and then unfold it
\textcolor{red}{to form a fixed angle $\theta$ at all creases}.
The resulting shape is called the Dragon curve with the unfolding 
angle $\theta$.
When $0\le\theta<90^{\circ}$, the corresponding Dragon curve has
a self-intersection.
When $\theta=180^{\circ}$, the corresponding Dragon curve
is a straight line, which has no self-intersection.
In this paper, we will show that any Dragon curve whose unfolding angle
is greater than $99.3438^{\circ}$ and less than $180^{\circ}$
has no self-intersection.
\end{abstract}

\section{Introduction}\label{introduction}

\begin{figure}
\caption{Folding-unfolding process for $k=1,2$}
\label{Paperfolding1}
\includegraphics[width=10cm]{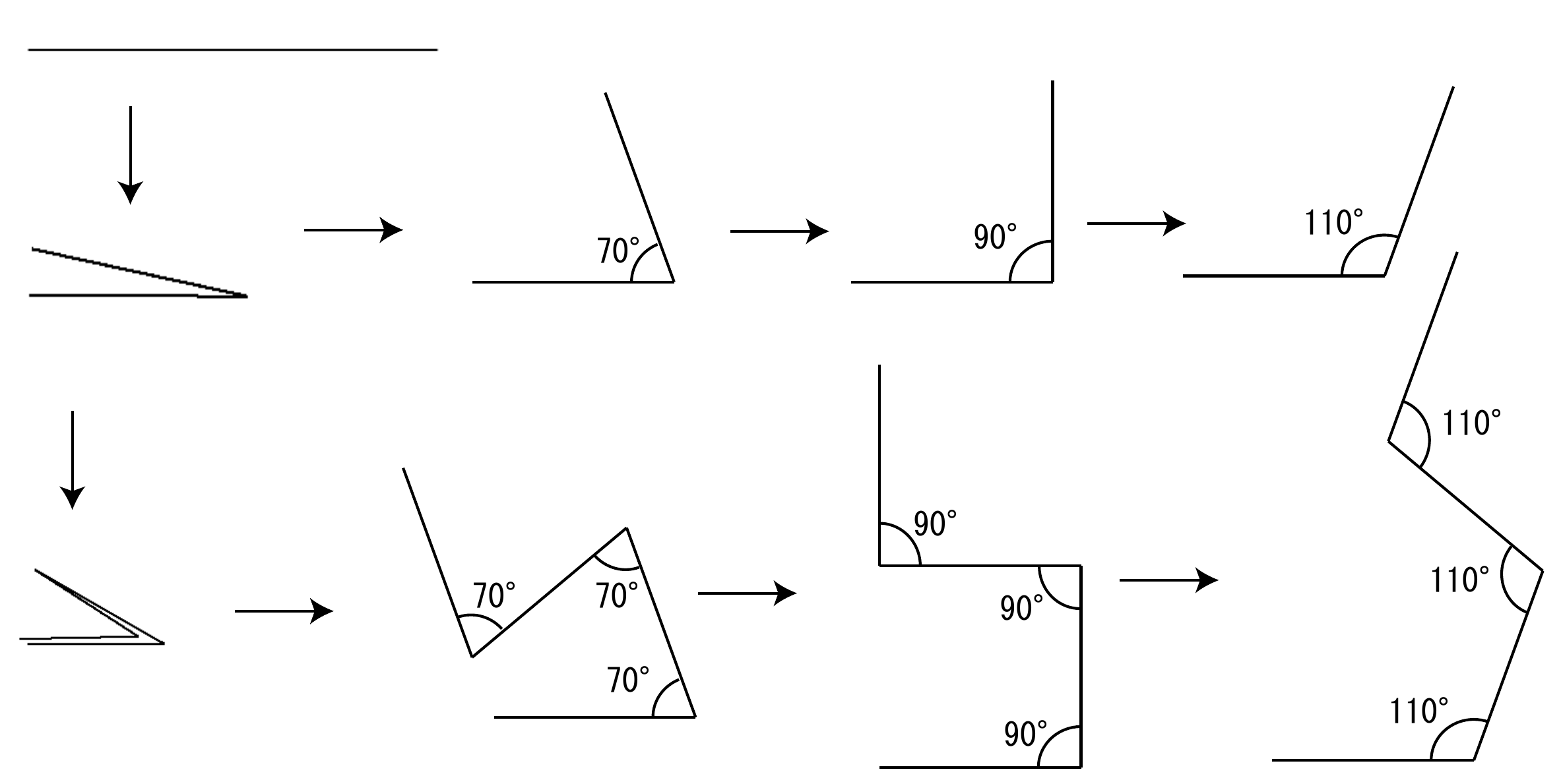}
\end{figure}

\begin{figure}
\begin{center}
\caption{The Heighway Dragon of order 8}
\label{Paperfolding2}
\includegraphics[clip, width=6.5cm]{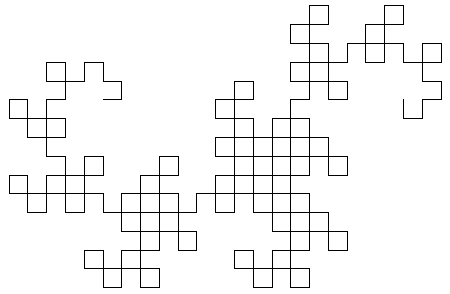}
\end{center}
\end{figure}

{\color{red}
Let us fold a strip of paper $k$ times in the same direction,
and then unfold it
\textcolor{red}{to form a fixed angle $\theta$ at all creases}.
We obtain the Dragon curve of order $k$ with the unfolding angle $\theta$ 
by this folding-unfolding process. 
Figure \ref{Paperfolding1} illustrates this process for $k=1$ and $k=2$.
In particular, the Dragon curve with the unfolding angle $\theta=90^{\circ}$ is known as the Heighway Dragon, see Figure \ref{Paperfolding2}.
Davis--Knuth \cite{DK} (see also Knuth \cite[pp.\,571--603]{Knu}
and its addendum \cite[pp.\,603--614]{Knu}) 
showed that 
the Heighway Dragon does not meet itself except when meeting between creases,
and eventually gives a plane-filling curve.
}

%







%
%

%








Tabachnikov \cite{Tab} gave a beautiful survey on Dragon curves 
and summerized remaining questions. 
Motivated by this article,
self-intersection properties of \textcolor{red}{Dragon curves} are revisited in
Allouche--Mend{\` e}s\;France--Skordev \cite{AFS},
with an effort to shed light on historically forgotten works. 
It is remarked in \cite[p.\,7, l.\,9--10]{AFS}
that
if $0\le \theta<90^{\circ}$, then the corresponding Dragon curve \textcolor{red}{of order 4} has
a self-intersection, because there are consecutive three turns to the same direction (right or left).
In the unpublished work of Albers \cite{Alb},
a less trivial fact is shown that if $90^{\circ}<\theta<95.126^{\circ}$
then the Dragon curve \textcolor{red}{of order 10} has a self-intersection.
\textcolor{red}{
Recently, \cite{Kamiya:21} showed that there is
a sequence $(\theta_n)_{n=1,2,\dots}$ converging to $\theta_0\approx 93.7912^{\circ}$
such that two creases of dragon curve with the unfolding angle $\theta_n$ exactly meet.}

In this paper, we say a Dragon curve is ``self-intersective''
if it meets itself at some point, i.e., it is not a Jordan curve.
Note that our ``self-intersective'' is slightly weaker than the one used
in \cite{AFS}. Even if it does not transversally \textcolor{red}{intersect} with itself
and \textcolor{red}{only meets} at \textcolor{red}{a crease}, we say it is self-intersective.
Hence, \textcolor{red}{the} Heighway Dragon is self-intersective in this paper and
the Dragon curve whose unfolding angle $\theta$ is in the range
$0^{\circ}\le \theta<95.126^{\circ}$ is self-intersective.

Under this terminology,
\textcolor{red}{\cite[Chapter 5.3]{AFS}} essentially asked a basic question;
{\it is there a constant $\eta$ with $95.126^{\circ}<\eta\le 180^{\circ}$
such that every Dragon curve whose unfolding angle $\theta$ is
in the range $\eta\le \theta\le 180^{\circ}$ is non\textcolor{red}{-}self-intersective?}
The statement is obvious when $\eta=180^{\circ}$;
the Dragon curve is a straight line. Apart from this trivial case,
there \textcolor{red}{seems to be} no rigorous proof
for non\textcolor{red}{-}self-intersectivity of any Dragon curve,
besides a lot of figures produced by computer experiments,
that ``clearly" suggest to be non\textcolor{red}{-}self-intersective.
Our paper is to fill this gap between computer-assisted intuition and mathematics.

Our method to address this question is to consider the renormalized fractal curve.
In the complex plane,
choose the similitude that sends
two endpoints of \textcolor{red}{a Dragon curve of order $k$} to $0$ and $1$. Then
we obtain an associated iterated function system (in short IFS).
Starting from a segment $[0,1]$ as an initial compact set,
\textcolor{red}{$k$ times} iteration of contractive maps for the IFS
gives the renormalized \textcolor{red}{Dragon curve of order $k$}.
Its Hausdorff limit as $k\to \infty$ is the attractor of the IFS.
To distinguish from the Dragon curve (broken line),
hereafter we call it the limit Dragon curve.
Let us define the limit Dragon curve and the corresponding
renormalized \textcolor{red}{Dragon curve of order $k$}.

\begin{defn}\label{def1}
Let $\xi$ be a fixed angle in $0\le \xi<\pi/3$.
Set $x$ and $\alpha$ to be $x=2\cos\xi$ and $\alpha=e^{-i\xi}/x$.
Define the functions $f_{1}$, $f_{2}$,
and $\psi$ on \textcolor{red}{the complex plane} by
$f_{1}(z)=\alpha z$, $f_{2}(z)=-\overline{\alpha} z +1$, and
\begin{equation}\label{defpsi}
\psi(z)=-\frac{\overline{\alpha}}{\alpha}(z-\alpha)+\alpha.
\end{equation}
It is easily verified that
\begin{equation}
\alpha+\overline{\alpha}=1, \label{alphakannkei}
\end{equation}
and
\begin{equation}
f_{2}=\psi\circ f_{1}. \label{f2nobunnkai}
\end{equation}
The limit Dragon curve ${\mathcal D}(\xi)$ is defined to be the attractor of
the ${\rm IFS}\{f_{1},f_{2}\}$.
The unfolding angle $\theta$
with respect to ${\mathcal D}(\xi)$ is $\theta=\pi -2\xi$.

Let $D_{0}$ be the segment $[0,1]$.
For any positive integer $k$,
the renormalized \textcolor{red}{Dragon curve of order $k$, $D_{k}$,}
is recursively defined by
\begin{equation}\label{broken}
D_{k}=f_{1}(D_{k-1}) \cup f_{2}(D_{k-1}).
\end{equation}
\end{defn}

Since the range of $|\alpha|$ is $1/2 \le |\alpha| < 1$ for $0\le \xi<\pi/3$,
the functions $f_{1}$ and $f_{2}$ are contractive and ${\mathcal D}(\xi)$
is well-defined as the attractor of the ${\rm IFS}$.
$D_{k}$ is the broken line consisting of $2^{k}$ segments,
and its Hausdorff limit as $k\to \infty$ is identical
\textcolor{red}{to} ${\mathcal D}(\xi)$.
Especially, ${\mathcal D}(\pi/4)$, whose unfolding angle is $\pi/2(=90^{\circ})$,
is the Hausdorff limit of the renormalization of \textcolor{red}{the} Heighway Dragon.
When $\xi=0$,
\textcolor{red}{it follows that $f_{1}(z)=z/2$, $f_{2}(z)=-z/2 +1$,
and $[0,1]=f_{1}([0,1]) \cup f_{2}([0,1])$.}
Hence the limit Dragon curve $\mathcal{D}(0)$ is the segment $[0,1]$,
\textcolor{red}{which has no self-intersection.}
To discuss the non\textcolor{red}{-}self-intersectivity of  $\mathcal{D}(\xi)$,
we restrict the range of $\xi$ to be $0<\xi<\pi/4$.

The non\textcolor{red}{-}self-intersectivity of $\mathcal{D}(\xi)$ is derived from
the following key proposition.

\begin{proposition}\label{Arc}
\textcolor{red}{Let $\xi$ be a fixed angle in $0< \xi<\pi/4$,}
$f_{1}$ and $f_{2}$ be the functions of Definition \ref{def1},
and ${\mathcal D}(\xi)$ be the limit Dragon curve
which is the attractor of the ${\rm IFS}\{f_{1},f_{2}\}$.
If there exists a non-empty open set $U$
\textcolor{red}{in the complex plane} which satisfies

{\rm (i)} $U\supset f_{1} (U) \cup f_{2}(U)$,

{\rm (ii)} $f_{1}(U)\cap f_{2}(U)=\emptyset$,

{\rm (iii)} $f_{1}(\overline{U})\cap f_{2}(\overline{U})=\{\alpha\}$,

\noindent where $\overline{U}$ is the closure of $U$ and $\alpha$ is the
same one as in Definition \ref{def1},
then ${\mathcal D}(\xi)$ is a simple arc, i.e., it is homeomorphic to the segment $[0,1]$.
\end{proposition}

The conditions (i) and (ii) are the open set condition
(see Moran \cite{Moran}, Hutchinson \cite{Hut}),
which implies the strong open set condition
(see Bandt--Graf \cite{BandtGraf:92}, Schief \cite{Schief:94}) in this case.
The inclusion (i) and the uniqueness of the attractor
give the inclusion ${\mathcal D}(\xi)\subset \overline{U}$.
The condition (iii) is necessary to show that ${\mathcal D}(\xi)$
is homeomorphic to $[0,1]$.
Proposition \ref{Arc} is not new. It often appears
in fractal geometry that under conditions (i), (ii), and (iii),
the attractor of an {\rm IFS} is homeomorphic to the segment $[0,1]$,
see e.g., Hata \cite{Hata:85} and Kigami \cite[p.\,15, Example 1.2.7]{Kigami}.
For the convenience of the reader, we give
a short proof of Proposition \ref{Arc} in Section \ref{homeo}.

We will prove that
there exists a constant $\xi_{0}$ with $0<\xi_{0}<\pi/4$ such that,
for any $\xi$ in the range $0<\xi<\xi_{0}$
and the corresponding limit Dragon curve $\mathcal{D}(\xi)$,
a certain non-empty open set $U_{\xi}$ satisfying
the conditions (i), (ii), and (iii) with $U=U_{\xi}$ can be chosen.
Such a non-empty open set $U_{\xi}$ is chosen
as the interior of a set $\mathcal{C}$ whose definition is
given in Section \ref{polygonC}.
In Section \ref{proofI}, it is proved that the condition (i) with $U=U_{\xi}$
is satisfied for any $\xi$ in $0<\xi<\pi/4$.
Sections \ref{sufficientII} and \ref{proofII}
are devoted to show that the condition (ii) with $U=U_{\xi}$
is satisfied for any $\xi$ in $0<\xi<\xi_{0}$.
The condition (iii) with $U=U_{\xi}$ is discussed in Section \ref{proofTh1}.
Then, by Proposition \ref{Arc}, we have the following.

\begin{theorem}\label{MainTheorem}
There exists a constant $\xi_{0}$, whose approximate value is
$\xi_{0}\approx 0.703858$, such that
the limit Dragon curve ${\mathcal D}(\xi)$ is a simple arc for any $\xi$
in \textcolor{red}{$0\le\xi<\xi_{0}$}.
\end{theorem}

\begin{figure}[H]
\begin{center}
\caption{The limit Dragon curve ${\mathcal D}(0.703858)$ is a simple arc.}
\label{ADragon}
\includegraphics[clip, width=6.5cm]{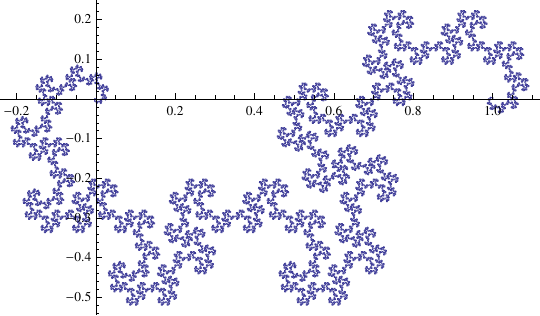} 	
\end{center}
\end{figure}

The approximate value of $\pi-2\xi_{0}$ is $1.73388$,
which is corresponding to the angle $99.3438^{\circ}$.
Hence any limit Dragon curve
whose unfolding angle is greater than $99.3438^{\circ}$
and less than \textcolor{red}{or equal to} $180^{\circ}$ is a \textcolor{red}{simple arc}.
Figure \ref{ADragon} is the limit Dragon curve ${\mathcal D}(0.703858)$
with the unfolding angle $99.3438^{\circ}$.

The technique of deriving Theorem \ref{MainTheorem} will be modified
in Section \ref{proofTh2} to show the \textcolor{red}{non-self-intersectivity}
of the renormalized \textcolor{red}{Dragon curve of order $k$.}
We have the following.

\begin{theorem}\label{Main2}
There exists a constant $\xi_{0}$,
whose approximate value is $\xi_{0}\approx 0.703858$,
such that the renormalized \textcolor{red}{Dragon curve of order $k$, $D_{k}$,}
has no self-intersection
for any positive integer $k$ and any $\xi$ in \textcolor{red}{$0\le\xi<\xi_{0}$}.
\end{theorem}


\section{Preliminaries}\label{preliminary}

Throughout the paper, ${\bf R}$ is the set of real numbers,
${\bf C}$ is the set of complex numbers, $i=\sqrt{-1}$,
$|z|$ is the usual norm of $z\in {\bf C}$,
and $\overline{z}$ is the complex conjugate of $z\in {\bf C}$.
For a set $A\subset {\bf C}$, $A^{\circ}$ is the interior of $A$,
and $\overline{A}$ is the closure of $A$.
For a point $p\in {\bf C}$ and $r>0$,
$B_{r}(p)$ (resp.\ $B_{r}^{\circ}(p)$)
is the closed (resp.\ open) ball with center $p$ and radius $r$.
For a point $p\in {\bf C}$ and a line (or segment) $L\subset {\bf C}$,
${\rm dist}(p,L)$ is defined by
${\rm dist}(p,L)=\min_{z\in L}\{|z-p|\}$.

For a finite word $a_{1}\cdots a_{k}$
\textcolor{red}{on the alphabet \{1, 2\}},
and for the functions $f_{1}$ and $f_{2}$ of Definition \ref{def1},
the function $f_{a_{1}\cdots a_{k}}(z)$
is defined by
$
f_{a_{1}\cdots a_{k}}(z)=
\textcolor{red}{(f_{a_{1}}\circ\cdots\circ f_{a_{k}})(z)}.
$
For a set $A\subset {\bf C}$, $f_{a_{1}\cdots a_{k}}(A)$ is the set defined by
$f_{a_{1}\cdots a_{k}}(A)=\{f_{a_{1}\cdots a_{k}}(z) : z\in A\}$.

For any one-sided infinite word $a_{1}a_{2}a_{3}\cdots$
\textcolor{red}{on the alphabet \{1, 2\}}
, the limit
$\lim_{k\to\infty}f_{a_{1}\cdots a_{k}}(z)
$
exists in ${\mathcal D}(\xi)$,
where the limit is independent of the choice of $z\in{\bf C}$.
In this paper we use the short notation $f_{a_{1}a_{2}a_{3}\cdots}$
for the above limit:
$f_{a_{1}a_{2}a_{3}\cdots}=\lim_{k\to\infty}f_{a_{1}\cdots a_{k}}(z)$.
For example, the points $0$, $1$, and $\alpha$
belong to ${\mathcal D}(\xi)$ with the expressions
$0=f_{111\cdots}$, $1=f_{2111\cdots}$, and $\alpha=f_{12111\cdots}=f_{22111\cdots}$.
A preperiodic word
$a_{1}\cdots a_{k}b_{1}\cdots b_{l}b_{1}\cdots b_{l}\cdots$
is denoted by $a_{1}\cdots a_{k}(b_{1}\cdots b_{l})^{\infty}$.
For example, $f_{111\cdots}=f_{(1)^{\infty}}$,
$f_{2111\cdots}=f_{2(1)^{\infty}}$, $f_{12111\cdots}=f_{12(1)^{\infty}}$,
and $f_{22111\cdots}=f_{22(1)^{\infty}}$.

The function $f_{a_{1}\cdots a_{k}}(z)$ is expressed as
\begin{equation}
f_{a_{1}\cdots a_{k}}(z)=\Big(\prod_{j=1}^{k}C_{a_{j}}\Big)\times
(z-f_{(a_{1}\cdots a_{k})^{\infty}})+f_{(a_{1}\cdots a_{k})^{\infty}}, \label{periodicfct}
\end{equation}
where $C_{a_{j}}$ is the coefficient of $z$ in $f_{a_{j}}(z)$,
because all of $f_{a_{j}}(z)$, $1\le j\le k$, are linear functions of $z$
and the fixed point of $f_{a_{1}\cdots a_{k}}(z)$ is $f_{(a_{1}\cdots a_{k})^{\infty}}$.
Substituting $z=0=f_{(1)^{\infty}}$ into \eqref{periodicfct}, we have
\begin{align}
& f_{(a_{1}\cdots a_{k})^{\infty}}
=\frac{f_{a_{1}\cdots a_{k}(1)^{\infty}}}{1-\prod_{j=1}^{k}C_{a_{j}}}, \label{periodicpt}\\
& f_{a_{1}\cdots a_{k}}(z)=\Big(\prod_{j=1}^{k}C_{a_{j}}\Big)\times z
+f_{a_{1}\cdots a_{k}(1)^{\infty}}. \label{linearexpress}
\end{align}
For example, since $f_{2(1)\infty}=1$ and $f_{12(1)\infty}=f_{22(1)\infty}=\alpha$,
we have, by \eqref{linearexpress},
\begin{align}
& f_{21}(z)=-|\alpha|^{2}z+1,\label{f21} \\
& f_{12}(z)=-|\alpha|^{2}z+\alpha,\label{f12} \\
& f_{22}(z)=(-\overline{\alpha})^{2}z+\alpha. \label{f22}
\end{align}

The point $f_{(2211)^{\infty}}$ plays an essential role in this paper.
Hereafter, $f_{(2211)^{\infty}}$ is denoted by $z_{0}$.
A geometric meaning of $z_{0}$ is explained in Remark \ref{rem4} below.
By \eqref{periodicpt}, $z_{0}$ is expressed as
\begin{equation}
z_{0}=f_{(2211)^{\infty}}=\frac{\alpha}{1-|\alpha|^{4}}. \label{z0}
\end{equation}

The range of the argument of a complex number $z$
is restricted to be $-\pi < \arg z\le \pi$.
The following properties for $\arg z$ are useful;
(i) $\arg(tz)=\arg{z}$ holds for any positive number $t$,
(ii) $0< \arg z<\pi$ is equivalent to $\Im z>0$,
(iii) if $0< \arg z<\pi$ holds, then
$0< \arg (z+t)<\pi$ holds for any real number $t$.

{\color{red}
For the triple $(a,b,c)$, where $a, b, c\in{\bf C}$ satisfy both conditions $a\neq b$ and $b\neq c$,
$\angle (a,b,c)$ is defined by
$\angle (a,b,c)=\arg((c-b)/(a-b))$.
}

\begin{defn}\label{def2}
{\color{red}

Let $a$ and $b$ be the distinct points in ${\bf C}$.
In the complex plane, directed lines {\rm (}half-lines{\rm )},
directed half-planes, closed line segments,
and polygonal convex sets are defined as follows.

\begin{enumerate}
\item[(i)] The directed line $L(a,b)$
is defined by the line through $a$ and $b$
which is directed from $a$ to $b$.
The directed half-line $HL(a,b)$ is defined by the half-line
through $a$ and $b$ which emanates from $a$.
\item[(ii)] The left {\rm (}resp.\ right{\rm )} half-plane $V^{+}(a,b)$
{\rm (}resp.\ $V^{-}(a,b)${\rm )}
with respect to $L(a,b)$ is defined by
$V^{+}(a,b)=\{z\in {\bf C} : 0< \angle (b,a,z)<\pi \}$
 {\rm (}resp.\
$V^{-}(a,b)=\{z\in {\bf C} : -\pi< \angle (b,a,z)<0 \}${\rm )}.
\textcolor{red}{The closed left {\rm (}resp.\ right{\rm )} half-plane $\overline{V^{+}(a,b)}$
{\rm (}resp.\ $\overline{V^{-}(a,b)}${\rm )} is defined to be the closure of
$V^{+}(a,b)$
{\rm (}resp.\ $V^{-}(a,b)${\rm )}.}
\item[(iii)] For the distinct points $v_{1}$ and $v_{2}$ in ${\bf C}$,
the closed line segment $[v_{1},v_{2}]$ is defined by
$[v_{1},v_{2}]=\{v_{1}+t(v_{2}-v_{1}) : 0\le t\le 1\}$.
We do not consider direction for closed line segments.
\item[(iv)] Let $(v_{n})_{n=1}^{N}$, $N\ge 3$, be a complex sequence such that
$0< \angle (v_{n-1},v_{n},v_{n+1})< \pi$ is satisfied for all $n$,
where $v_{0}$ and $v_{N+1}$ is defined to be $v_{0}=v_{N}$ and $v_{N+1}=v_{1}$.
Then the convex set
$
\mathcal{P}((v_{n})_{n=1}^{N})
$
associated with $(v_{n})_{n=1}^{N}$ is defined by
$$
\mathcal{P}((v_{n})_{n=1}^{N})=\bigcap_{n=1}^{N}
\overline{V^{+}(v_{n+1},v_{n})}.
$$
\end{enumerate}
}
\end{defn}

\begin{remark}\label{rem1}
For any linear function $f(z)$ on ${\bf C}$,
$\angle (f(a),f(b),f(c))=\angle (a,b,c)$ holds because
$(f(c)-f(b))/(f(a)-f(b))=(c-b)/(a-b)$.
\end{remark}

\begin{remark}\label{rem2}
The definition of $V^{+}(a,b)$ is independent of the choice of $a,b\in L(a,b)$.
In fact, for $a',b'\in L(a,b)$ whose expressions are
$a'=a+t_{a'}(b-a)$, $b'=a+t_{b'}(b-a)$ with $t_{a'}<t_{b'}$,
the property {\rm (i)} for $\arg z$ gives
$$
\angle(b',a',z)=\arg\Big(\frac{1}{t_{b'}-t_{a'}}(\frac{z-a}{b-a}-t_{a'})\Big)
=\arg\Big(\frac{z-a}{b-a}-t_{a'}\Big),
$$
then the property {\rm (iii)} for $\arg z$
and the definition of $V^{+}(a,b)$ give $V^{+}(a,b)=V^{+}(a',b')$.
\end{remark}

\begin{remark}\label{rem3}
To show inclusion such as
$
\mathcal{P}((u_{m})_{m=1}^{M})\subset\mathcal{P}((v_{n})_{n=1}^{N})
$
or
$[u_{1},u_{2}]\subset\mathcal{P}((v_{n})_{n=1}^{N})$,
we only need to show
$u_{m}\in \mathcal{P}((v_{n})_{n=1}^{N})$ for all $m$
by the convexity property of $\mathcal{P}((v_{n})_{n=1}^{N})$.
\end{remark}

\begin{remark}\label{rem4}
The point $z_{0}$ has been already introduced as $z_{0}=f_{(2211)^{\infty}}$.
Let us consider the sequence
$(\alpha,f_{2211}(\alpha),f_{22112211}(\alpha),
f_{221122112211}(\alpha),\ldots)$.
Then $z_{0}$ is the limit point for the sequence.
$f_{2211}(\alpha)$ is an \textcolor{red}{external division point of} $[0,\alpha]$.
This is verified by calculating the angle
$\angle (0,\alpha,f_{2211}(\alpha))$ as
\begin{align*}
\angle (0,\alpha,f_{2211}(\alpha))	
= \arg \frac{f_{2211}(\alpha)-\alpha}{0-\alpha}
= \arg (-1\times|\alpha|^{4})
= \arg (-1)=\pi,
\end{align*}
where the second equality comes from \eqref{f22}
and the third equality comes from the property {\rm (i)} for $\arg z$.
It is similarly verified that all points in the sequence
and the limit point $z_{0}$ are on the directed half-line $HL(0,\alpha)$,
and the distance between $0$ and $z_{0}$ is greater than
that between $0$ and any point in the sequence.

The point $z_{0}$ has been introduced by the third author in \cite{Wen:22} as
the explicit extremal point of the convex hull of
the limit Dragon curve ${\mathcal D}(\xi)$.
\end{remark}

\section{Proof of Proposition \ref{Arc}}\label{homeo}

We prove Proposition \ref{Arc} following the
notations in \cite[pp.\,12--15]{Kigami}.
Let $\Sigma$ be the set of one-sided infinite words
\textcolor{red}{on the alphabet \{1, 2\}}
, that is,
$\Sigma=\{\omega : \omega=\omega_{1}\omega_{2}\omega_{3}\cdots, \omega_{j}\in\{1,2\}\}$.
Let $\delta_{r}$, $0<r<1$, be the metric on $\Sigma$ defined by
$\delta_{r}=r^{s(\omega,\tau)}$, where
$\omega=\omega_{1}\omega_{2}\omega_{3}\cdots, \tau=\tau_{1}\tau_{2}\tau_{3}\cdots
\in\Sigma$ and \textcolor{red}{$s(\omega,\tau)=0$ if $\omega=\tau$} and $s(\omega,\tau)=\min\{m:\omega_{m}\neq\tau_{m}\}-1$ \textcolor{red}{otherwise}.
Then it is known that $(\Sigma,\delta_{r})$ is a compact metric space
(see Theorem 1.2.2 of \cite{Kigami}).

For the ${\rm IFS}\{f_{1},f_{2}\}$ and its attractor ${\mathcal D}(\xi)$
of Definition \ref{def1}, and for any
$\omega=\omega_{1}\omega_{2}\omega_{3}\cdots\in\Sigma$,
$
\bigcap_{m=1}^{\infty}f_{\omega_{1}\cdots\omega_{m}}({\mathcal D}(\xi))
$
consists of only one point, \textcolor{red}{which is written $f_{\omega}$} in this paper.
Define the map $\pi_{\xi} : \Sigma\to{\mathcal D}(\xi)$ by
$\pi_{\xi}(\omega)=f_{\omega}$.
Then it is known that the map $\pi_{\xi}$ is continuous and surjective,
where the metric on ${\mathcal D}(\xi)$ comes from the usual metric on ${\bf C}$
(see Theorem 1.2.3 of \cite{Kigami}).

By Proposition 1.2.5 of \cite{Kigami},
the equality $\pi_{\xi}(\omega)=\pi_{\xi}(\tau)$, $\omega\neq\tau$, holds
if and only if
$\pi_{\xi}(\sigma^{s(\omega,\tau)}\omega)=\pi_{\xi}(\sigma^{s(\omega,\tau)}\tau)$
holds, where $\sigma$ is the shift map.
The inclusion (i) of Proposition \ref{Arc}
and the uniqueness of the attractor
give the inclusion ${\mathcal D}(\xi)\subset \overline{U}$.
Then, by this inclusion and (iii) of Proposition \ref{Arc},
$
\pi_{\xi}(\sigma^{s(\omega,\tau)}\omega)=\pi_{\xi}(\sigma^{s(\omega,\tau)}\tau)
\in
f_{1}({\mathcal D}(\xi))\cap f_{2}({\mathcal D}(\xi))
\subset f_{1}(\overline{U})\cap f_{2}(\overline{U})=\{\alpha\}
$
holds. Hence, since $\alpha$ is expressed as
$\alpha=f_{12(1)^{\infty}}=f_{22(1)^{\infty}}$,
the equality $\pi_{\xi}(\omega)=\pi_{\xi}(\tau)$, $\omega\neq\tau$, holds
if and only if there exists a finite word $w$ such that
$\{\omega,\tau\}=\{w12(1)^{\infty},w22(1)^{\infty}\}$.
The natural identification $w12(1)^{\infty}\sim w22(1)^{\infty}$
gives the quotient space $\Sigma /\sim$.
Since the space $\Sigma$ is compact,
the quotient space $\Sigma / \sim$ is compact.
Let $p: \Sigma\to\Sigma / \sim$ be the natural projection, and
define the map
$\overline{\pi_{\xi}}: \Sigma / \sim \to {\mathcal D}(\xi)$
by $\pi_{\xi}=\overline{\pi_{\xi}}\circ p$.
Then the map $\overline{\pi_{\xi}}$ is bijective.
Moreover, $\overline{\pi_{\xi}}$ is a closed map because of the
compactness of $\Sigma / \sim$.
Thus, $\overline{\pi_{\xi}}$ is a continuous, bijective,
and closed map, and consequently, is a homeomorphism.

{\color{red}
Define the functions $\varphi_{1}$ and $\varphi_{2}$
of a real variable $y$ by $\varphi_{1}(z)=y/2$ and $\varphi_{2}(z)=-y/2 +1$.
The attractor of the ${\rm IFS}\{\varphi_{1},\varphi_{2}\}$ is the segment $[0,1]$
in ${\bf R}$, and the open segment $(0,1)$ in ${\bf R}$ satisfies the conditions
(i) $(0,1)\supset \varphi_{1} ((0,1)) \cup \varphi_{2}((0,1))$,
(ii) $\varphi_{1} ((0,1)) \cap \varphi_{2}((0,1))=\emptyset$, and
(iii) $\varphi_{1}([0,1])\cap \varphi_{2}([0,1])=\{ 1/2 \}$.
Since a similar discussion works on the map
$\overline{\pi_{0}}: \Sigma / \sim \to [0,1]$,
$\overline{\pi_{0}}$ is a homeomorphism.
Thus $\overline{\pi_{\xi}}\circ\overline{\pi_{0}}^{-1}$ is a
homeomorphism between $[0,1]$ and $\mathcal{D}(\xi)$.
}
This completes the proof of Proposition \ref{Arc}.


\section{An infinite polygon with respect to
${\mathcal D}(\xi)$}\label{polygonC}

In this section, an infinite polygon $\mathcal{C}$ is introduced
(see Definition \ref{def4} below).  The interior of $\mathcal{C}$
is chosen as $U_{\xi}$ in Section 1.
The point $z_{0}=f_{(2211)^{\infty}}$ is one of the vertices of $\mathcal{C}$.
To prove Theorem \ref{MainTheorem}, we will verify that
the conditions $\mathcal{C}\supset f_{1} (\mathcal{C})\cup f_{2}(\mathcal{C})$,
$f_{1}(\mathcal{C})\cap f_{2}(\mathcal{C})=\emptyset$, and
$f_{1}(\overline{\mathcal{C}})\cap f_{2}(\overline{\mathcal{C}})=\{\alpha\}$
are satisfied.
The choice $\mathcal{C}$ will not be the best possible,
but we have to take a balance between the computational complexity
and the accuracy of the resulting bound.

\begin{figure}[H]
\begin{center}
\caption{The sets $A_{1}$ and $\widetilde{A_{1}}$}
\label{A1}
\includegraphics[clip, width=8.cm]{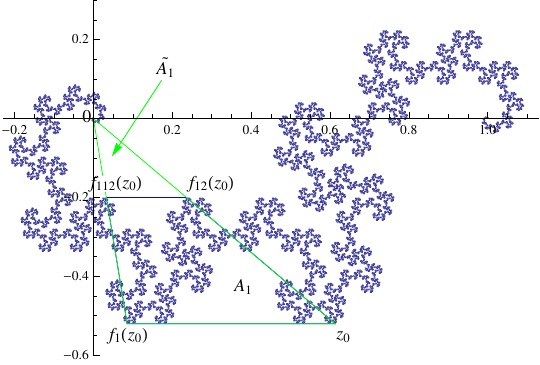}
\end{center}
\end{figure}

\begin{defn}\label{def3}
Let $z_{0}$ be the point of \eqref{z0},
and $f_{1}$ and $f_{2}$ be the functions of Definition \ref{def1}.
For an integer $m$, define the function $f_{(1)^{m}}(z)$ by
$f_{(1)^{m}}(z)=f_{\scriptsize{\underbrace{1\cdots 1}_{m}}}(z)$
if $m$ is positive,
$f_{(1)^{m}}(z)=f_{(1)^{-m}}^{-1}(z)$
if $m$ is negative, and $f_{(1)^{0}}(z)=z$.
The sets $A_{1}$ and $\widetilde{A_{1}}$ are defined by
$A_{1}=\mathcal{P}(z_{0},f_{1}(z_{0}),f_{112}(z_{0}),f_{12}(z_{0}))$
and $\widetilde{A_{1}}=\mathcal{P}(0,z_{0},f_{1}(z_{0}))$ respectively,
see Figure \ref{A1}.
The sets $A_{m}$ and $\widetilde{A_{m}}$ are defined by
$A_{m}=f_{(1)^{m-1}}(A_{1})$ and
$\widetilde{A_{m}}=f_{(1)^{m-1}}(\widetilde{A_{1}})$ respectively.
\end{defn}

\vspace{-9mm}

\begin{remark}\label{rem5}
We have already restricted the range of $\xi$ to be $0<\xi<\pi/4$.
This range of $\xi$ gives
the inequalities $1/2<|\alpha|<1/\sqrt{2}$ and
\begin{equation}\label{lower1}
	1-|\alpha|^{2}-|\alpha|^{4}>\frac{1}{4}.
\end{equation}
We can calculate angles related with $A_{1}$ and $\widetilde{A_{1}}$.
For example, by \eqref{f12}, \eqref{z0}, \eqref{lower1}, and
the property {\rm (i)} for $\arg z$,
\begin{align*}
\angle (0,f_{12}(z_{0}),z_{0})
= \arg \frac{z_{0}-f_{12}(z_{0})}{0-f_{12}(z_{0})}
= \arg \Big(-1\times
\frac{|\alpha|^{2}+|\alpha|^{4}}{1-|\alpha|^{2}-|\alpha|^{4}}\Big)
= \arg (-1)=\pi,
\end{align*}
then, by Remark \ref{rem1},
$\angle (0,f_{112}(z_{0}),f_{1}(z_{0}))= \angle (0,f_{12}(z_{0}),z_{0})=\pi$.
Hence $f_{12}(z_{0})$ {\rm (}resp. $f_{112}(z_{0})${\rm )}
\textcolor{red}{is an internal division point of} $[0,z_{0}]$ {\rm (}resp. $[0,f_{1}(z_{0})]${\rm )}.
We similarly have
$\angle (z_{0},f_{1}(z_{0}),f_{112}(z_{0}))=\theta$,
$\angle (f_{1}(z_{0}),f_{112}(z_{0}),f_{12}(z_{0}))=2\xi$,
$\angle (f_{112}(z_{0}),f_{12}(z_{0}),z_{0})=\theta +\xi$,
and $\angle (f_{12}(z_{0}),z_{0},f_{1}(z_{0}))=\xi$,
which are all greater than $0$ and less than $\pi$.
Thus $A_{1}$ is well-defined as a convex set
in the sense of Definition \ref{def2}.
The angles for $\widetilde{A_{1}}$ are
$\angle (0,z_{0},f_{1}(z_{0}))=\xi$,
$\angle (z_{0},f_{1}(z_{0}),0)=\theta$,
and $\angle (f_{1}(z_{0}),0,z_{0})=\xi$.
\textcolor{red}{The inclusion $A_{m}\subset \widetilde{A_{m}}$} holds for any $m$.
\end{remark}

\begin{figure}[H]
\begin{center}
\caption{The infinite polygon $\mathcal{C}$}
\label{Cover}
\includegraphics[clip, width=8.cm]{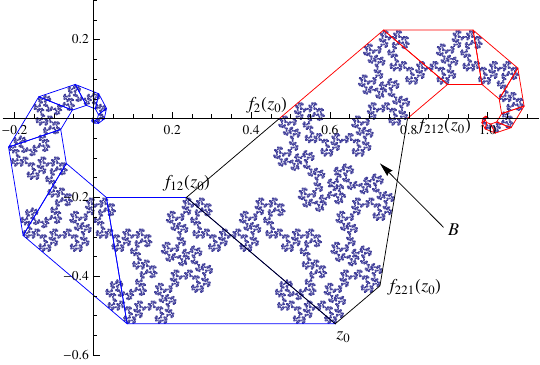}
\end{center}
\end{figure}

\begin{defn}\label{def4}
Let $z_{0}$ be the point of \eqref{z0},
$f_{1}$ and $f_{2}$ be the functions of Definition \ref{def1},
and $A_{m}$ be the set of Definition \ref{def3}.
Define the set $B$ by
\begin{align*}
B =A_{0}\cap \overline{V^{+}(f_{221}(z_{0}),f_{212}(z_{0}))} =\mathcal{P}(z_{0},f_{12}(z_{0}),f_{2}(z_{0}),f_{212}(z_{0}),f_{221}(z_{0})),
\end{align*}
see Figure \ref{Cover}.
Then the set $\mathcal{C}$ is defined by
$$
\mathcal{C}=\Big(\bigcup_{n=1}^{\infty}A_{n}\Big)
\cup B \cup \Big(\bigcup_{n=1}^{\infty}f_{2}(A_{n})\Big).
$$
\end{defn}

\section{Proof of $\mathcal{C}\supset f_{1} (\mathcal{C})\cup f_{2}(\mathcal{C})$
for $0<\xi<\pi/4$}\label{proofI}

In this section, we will prove that
$\mathcal{C}\supset f_{1} (\mathcal{C})\cup f_{2}(\mathcal{C})$
is satisfied for any $\xi$ in $0<\xi<\pi/4$.
Firstly, let us introduce simple polygons $S$, $S'$, and $S''$
which are including $\bigcup_{n=1}^{\infty}A_{n}$.

\begin{figure}[H]
\begin{center}
\caption{The points $p_{1}$, $p_{2}$, $p_{3}$, and $q$}
\label{Hishigata}
\includegraphics[clip, width=9.cm]{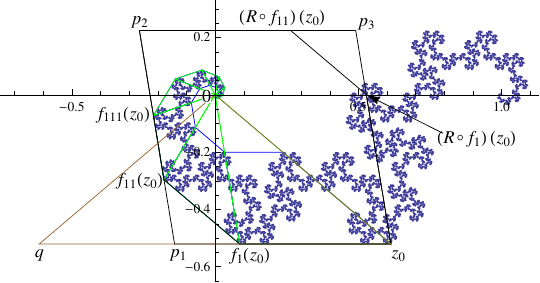}
\end{center}
\end{figure}

\begin{defn}\label{def5}
The points $p_{1},p_{2},p_{3}\in{\bf C}$ in Figure \ref{Hishigata}
are chosen as follows:

\textcolor{red}{{\rm (i)}}  $p_{1}$ is the intersection
of the lines $L(f_{1}(z_{0}),z_{0})$ and $L(f_{111}(z_{0}),f_{11}(z_{0}))$,

\textcolor{red}{{\rm (ii)}}  $p_{2}$ is the intersection of the lines
$L(f_{111}(z_{0}),f_{11}(z_{0}))$ and $L(z_{0},0)$,

\textcolor{red}{{\rm (iii)}}  $p_{3}$ satisfies both \textcolor{red}{conditions}
$\angle (z_{0},p_{2},p_{3})=\xi$ and $\angle (p_{3},z_{0},p_{2})=\xi$.

The set $S$ is defined by $S=\mathcal{P}(z_{0},p_{1},p_{2},p_{3})$
whose angles are
$\angle (z_{0},p_{1},p_{2})=\angle (p_{2},p_{3},z_{0})=\theta$ and
$\angle (p_{1},p_{2},p_{3})=\angle (p_{3},z_{0},p_{1})=2\xi$.
The set $S'$ is defined by
$S'=\mathcal{P}(z_{0},f_{1}(z_{0}),f_{11}(z_{0}),p_{2},p_{3})$.
Let $R(z)$ be the reflection with respect to $L(z_{0},p_{2})$
so that $R(p_{1})=p_{3}$ is satisfied, that is,
$R(z)=(\alpha/\overline{\alpha})\overline{z}$.
The set $S''$ is defined by
$S''=\mathcal{P}(z_{0},f_{1}(z_{0}),f_{11}(z_{0}),p_{2},
\textcolor{red}{(} R\circ f_{11}\textcolor{red}{)}(z_{0}),
\textcolor{red}{(} R\circ f_{1}\textcolor{red}{)}(z_{0}))$.

Choose the point $q\in{\bf C}$ in Figure \ref{Hishigata} so that
both \textcolor{red}{conditions}
$\angle (0,z_{0},q)=\xi$ and $\angle (z_{0},q,0)=\xi$ are satisfied.
The set $W$ is defined by $W=\mathcal{P}(0,z_{0},q)$
whose angles are
$\angle (0,z_{0},q)=\angle (z_{0},q,0)=\xi$ and $\angle (q,0,z_{0})=\theta$.
\end{defn}

\begin{lemma}\label{lemma1}
Under the same \textcolor{red}{notation} as in Definition \ref{def5},
we have the \textcolor{red}{following};

{\rm (i)} $q=-\frac{\overline{\alpha}}{\alpha}z_{0}
=z_{0}-\frac{1}{1-|\alpha|^{4}}$,

{\rm (ii)} $p_{1}=(\alpha-\overline{\alpha}|\alpha|^{2})z_{0}
=z_{0}-\frac{|\alpha|^{2}}{1-|\alpha|^{2}}$,

{\rm (iii)} $p_{2}=-|\alpha|^{2}z_{0}=f_{11}(q)$,

{\rm (iv)} $p_{3}=(\overline{\alpha}-\alpha|\alpha|^{2})z_{0}
=p_{2}+\frac{|\alpha|^{2}}{1-|\alpha|^{2}}$,

{\rm (v)} $f_{11}(W)\subset S''$.
\end{lemma}

\noindent{\it Proof.}
Define the function $\varphi_{1}$ by
$\varphi_{1}(z)=\alpha^{-1}(z-z_{0})+z_{0}$.
Then (i) is obtained by the relation $q=\varphi_{1}(0)$ and \eqref{alphakannkei}.
Define the function $\varphi_{2}$ by
$\varphi_{2}(z)=\overline{\alpha}(z-f_{1}(z_{0}))+f_{1}(z_{0})$.
Then (ii) is obtained by the relation $p_{1}=\varphi_{2}(f_{11}(z_{0}))$
and \eqref{alphakannkei}.
Define the function $\varphi_{3}$ by
$\varphi_{3}(z)=\overline{\alpha}^{-1}(z-z_{0})+z_{0}$.
Then (iii) is obtained by the relation $p_{2}=\varphi_{3}(p_{1})$,
the second equality of (ii), and \eqref{z0}.
Define the function $\varphi_{4}$ by
$\varphi_{4}(z)=\alpha(z-z_{0})+z_{0}$.
Then (iv) is obtained by the relation $p_{3}=\varphi_{4}(p_{2})$,
the first equality of (iii), and \eqref{alphakannkei}.
For (v), by (iii), it follows that
$f_{11}(W)=\mathcal{P}(0,f_{11}(z_{0}),p_{2})$,
for which all vertices $0$, $f_{11}(z_{0})$, and $p_{2}$
belong to $S''$. Hence (v) is obtained. \qed

\begin{remark}\label{rem6}
By the expressions of the points in Lemma \ref{lemma1},
it is verified that
\textcolor{red}{
\begin{align*}
& \angle (z_{0},f_{1}(z_{0}),p_{1})
=\angle (p_{2},f_{11}(z_{0}),p_{1})
=\angle (p_{2},(R\circ f_{11})(z_{0}),p_{3})\\
&= \angle (z_{0},(R\circ f_{1})(z_{0}),p_{3})
=\arg (-1\times |\alpha|^{2})
= \arg (-1)=\pi.
\end{align*}
}
\noindent Hence $S'' \subset S' \subset S$ holds.
\end{remark}

\begin{lemma}\label{lemma2}
Let $\xi$ be a fixed angle in $0<\xi<\pi/4$.
Choose a positive integer $N\ge 3$ so that
the inequality $\pi/(N+2)\le \xi< \pi/(N+1)$ is satisfied,
which is equivalent to the inequality $(N-1)\xi< \theta\le N\xi$.
Under the same \textcolor{red}{notation} as in Definition \ref{def5},
we have the \textcolor{red}{following};

{\rm (i)}
$f_{(1)^{n}}(z_{0})\in W\setminus [0,q]$
for all $n$ with $1\le n\le N-1$,

{\rm (ii)} $\angle (f_{(1)^{N}}(z_{0}),0,q)\ge 0$ and
$\angle (-1,0,f_{(1)^{N}}(z_{0}))>0$,

{\rm (iii)} $f_{(1)^{n}}(z_{0})\in B_{|f_{(1)^{N}}(z_{0})|}(0)$
for all $n$ with $n\ge N$,

{\rm (iv)} $|f_{(1)^{N}}(z_{0})|
<{\rm dist}(0,L(f_{1}(z_{0}),z_{0}))$.
\end{lemma}

\noindent{\it Proof.}
(i) From Remark \ref{rem6} and Lemma \ref{lemma1} (i) (ii),
it follows that
$f_{1}(z_{0})\in [p_{1},z_{0}]\subset [q,z_{0}]\setminus \{q\}$.
Hence $f_{1}(z_{0})\in W\setminus [0,q]$.
Let $h_{n}\in{\bf C}$, $1\le n\le N-1$, be the point such that
$
\{ h_{n} \}=
L(f_{(1)^{n}}(z_{0}), f_{(1)^{n-1}}(z_{0}))\cap [0,q],
$
where $h_{1}$ is equal to the point $q$.

Firstly, let us consider the triangle $\mathcal{P}(0,f_{1}(z_{0}),h_{1})$.
The angle $\angle (h_{1},0,f_{1}(z_{0}))$ is divided into
the angles $\angle (h_{1},0,f_{11}(z_{0}))=\theta-2\xi$
and $\angle (f_{11}(z_{0}),0,f_{1}(z_{0}))=\xi$.
The angle $\angle (0,f_{1}(z_{0}),h_{1})$ is divided into
the angles $\angle (0,f_{1}(z_{0}),f_{11}(z_{0}))=\xi$
and $\angle (f_{11}(z_{0}),f_{1}(z_{0}),h_{1})=\xi$.
Hence $f_{11}(z_{0})$ belongs to the interior of
the triangle $\mathcal{P}(0,f_{1}(z_{0}),h_{1})$.
On the other hand, since all vertices $0$,
$f_{1}(z_{0})$, and $h_{1}$ belong to $W$,
$\mathcal{P}(0,f_{1}(z_{0}),h_{1})$ is contained in $W$.
Thus $f_{11}(z_{0})\in \mathcal{P}^{\circ}(0,f_{1}(z_{0}),h_{1})
\subset W\setminus [0,q]$.

Next, let us consider the triangle $\mathcal{P}(0,f_{11}(z_{0}),h_{2})$
in the case $N\ge 4$.
The angle $\angle (h_{2},0,f_{11}(z_{0}))$ is divided into
the angles $\angle (h_{2},0,f_{111}(z_{0}))=\theta-3\xi$ and
$\angle (f_{111}(z_{0}),0,f_{11}(z_{0}))=\xi$.
The angle $\angle (0,f_{11}(z_{0}),h_{2})$) is divided into
the angles $\angle (0,f_{11}(z_{0}),f_{111}(z_{0}))=\xi$
and $\angle (f_{111}(z_{0}),f_{11}(z_{0}),h_{2})=\xi$.
Hence $f_{111}(z_{0})$ belongs to the interior of
the triangle $\mathcal{P}(0,f_{11}(z_{0}),h_{2})$.
On the other hand, since all vertices $0$,
$f_{11}(z_{0})$, and $h_{2}$ belong to $W$,
$\mathcal{P}(0,f_{11}(z_{0}),h_{2})$ is contained in $W$.
Thus $f_{111}(z_{0})\in \mathcal{P}^{\circ}(0,f_{11}(z_{0}),h_{2})
\subset W\setminus [0,q]$.

We can repeat this argument until $n=N-1$ because of $\theta-(N-1)\xi>0$.
Thus $f_{(1)^{N-1}}(z_{0})\in W\setminus [0,q]$.

(ii)
We have
$\angle (f_{(1)^{N}}(z_{0}),0,q)
=\angle (f_{(1)^{N}}(z_{0}),0,z_{0})-\angle (q,0,z_{0})
=N\xi -\theta\ge 0$ and
$
\angle (-1,0,f_{(1)^{N}}(z_{0}))
= \arg (-\alpha ^{N+1})
= \arg e^{i(\pi- (N+1)\xi)}=\pi - (N+1)\xi
= \theta+2\xi-(N+1)\xi=\theta-(N-1)\xi>0
$.

(iii) This is equivalent to
$|f_{(1)^{n-N}}(z_{0})|\le |z_{0}|$,
which follows from $|\alpha|<1/\sqrt{2}$.

(iv) We have
$|f_{(1)^{N}}(z_{0})|=|z_{0}|/(2\cos\xi)^{N}$ and
${\rm dist}(0,L(f_{1}(z_{0}),z_{0}))=|z_{0}|\sin \xi$.
Hence (iv) is equivalent to $1<(2\cos\xi)^{N}\sin \xi$.
The inequality $\pi/(N+2)\le \xi< \pi/(N+1)$ and $N\ge 3$ give
$2\cos\xi>2\cos(\pi/(N+1))\ge \sqrt{2}$ and
$\sin\xi\ge \sin(\pi/(N+2))$.
On the other hand, by convexity of $\sin x$ for $0<x<\pi/2$,
it follows that $\sin x>((\sin (\pi/5))/(\pi/5))x$ for $0<x<\pi/5$.
Hence we have
\begin{align*}
(2\cos\xi)^{N}\sin \xi
& >(\sqrt{2})^{N}\sin\frac{\pi}{N+2}\\
& > \min_{N\ge 3}\Big\{
(\sqrt{2})^{N}\frac{\sin \frac{\pi}{5}}{\frac{\pi}{5}}\frac{\pi}{N+2}\Big\}\\
& = \min_{M\ge 5}\Big\{\frac{5}{2}\frac{(\sqrt{2})^{M}}{M}\sin \frac{\pi}{5}\Big\}\\
& =2\sqrt{2}\sin \frac{\pi}{5}\approx 1.66251,
\end{align*}
where we use that $(\sqrt{2})^{M}/M$, $M\ge 5$, takes the minimum at $M=5$.
Thus $1<(2\cos\xi)^{N}\sin \xi$ is satisfied, and (iv) is obtained. \qed

\begin{lemma}\label{lemma3}
Let $\xi$ be a fixed angle in $0<\xi<\pi/4$,
$A_{n}$ and $\widetilde{A_{n}}$, $n\ge 1$, be the sets of Definition \ref{def3},
and $S$, $S'$, $S''$ be the sets of Definition \ref{def5}. Then we have
$$
\bigcup_{n=1}^{\infty}A_{n}
\subset \bigcup_{n=1}^{\infty}\widetilde{A_{n}}
\subset S'' \subset S' \subset S.
$$
\end{lemma}

\noindent{\it Proof.}
By Remark \ref{rem6}, it \textcolor{red}{suffices} to prove
$\bigcup_{n=1}^{\infty}\widetilde{A_{n}}\subset S''$.
Since all vertices $0$, $z_{0}$, and $f_{1}(z_{0})$ of $\widetilde{A_{1}}$
(resp.\ $0$, $f_{1}(z_{0})$, and $f_{11}(z_{0})$ of $\widetilde{A_{2}}$)
belong to $S''$, $\widetilde{A_{1}}$ (resp.\ $\widetilde{A_{2}}$)
is contained in $S''$.
We shall prove $\widetilde{A_{n}}\subset S''$ for $n\ge 3$.
Apply $f_{11}$ to Lemma \ref{lemma2} (i) (iii) (iv) to get
\vspace{3mm}

{\rm (i)'}\ \ $f_{(1)^{n}}(z_{0})\in f_{11}(W\setminus [0,q])
\subset f_{11}(W)$
for all $n$ with $3\le n\le N+1$,

{\rm (iii)'}\ \ $f_{(1)^{n}}(z_{0})\in
B_{|f_{(1)^{N+2}}(z_{0})|}(0)$
for all $n$ with $n\ge N+2$,

{\rm (iv)'}\ \ $|f_{(1)^{N+2}}(z_{0})|
<{\rm dist}(0,L(f_{111}(z_{0}),f_{11}(z_{0})))
={\rm dist}(0,[p_{2},f_{11}(z_{0})])$.
\vspace{3mm}

The triangle $\widetilde{A_{n}}$ is expressed as
$\widetilde{A_{n}}
=\mathcal{P}(0,f_{(1)^{n-1}}(z_{0}),f_{(1)^{n}}(z_{0}))$.
If $3\le n\le N+1$, then all vertices of $\widetilde{A_{n}}$ belong
to $f_{11}(W)$ by (i)' and $f_{(1)^{2}}(z_{0})\in f_{11}(W)$.
Hence $\widetilde{A_{n}}\subset f_{11}(W)$
for $3\le n\le N+1$. This and Lemma \ref{lemma1} (v) give
\begin{equation}\label{inclusionL}
\widetilde{A_{n}}\subset S''\quad \mbox{for $3\le n\le N+1$.}
\end{equation}

Let us consider $\widetilde{A_{n}}$ for $n\ge N+3$.
By Pythagorean theorem,
\begin{equation}\label{distS2L}
{\rm dist}(0,[p_{2},f_{11}(z_{0})])
< {\rm dist}(0,[f_{11}(z_{0}),f_{1}(z_{0})])
< {\rm dist}(0,[f_{1}(z_{0}),z_{0}]).
\end{equation}
By applying the reflection $R$ of Definition \ref{def5} to \eqref{distS2L},
\begin{align}
& {\rm dist}(0,[p_{2},f_{11}(z_{0})])
= {\rm dist}(0,[p_{2},
\textcolor{red}{(} R\circ f_{11}\textcolor{red}{)}(z_{0})]) \nonumber \\
& \quad < {\rm dist}(0,[\textcolor{red}{(}R\circ f_{11}\textcolor{red}{)} (z_{0}),
\textcolor{red}{(} R\circ f_{1}\textcolor{red}{)} (z_{0})])
< {\rm dist}(0,[\textcolor{red}{(} R\circ f_{1}\textcolor{red}{)}(z_{0}),z_{0}]). \label{distS2R}
\end{align}
\eqref{distS2L} and \eqref{distS2R} show that
the distance between $0$ and the boundary of $S''$ is equal to
${\rm dist}(0,[p_{2},f_{11}(z_{0})])$.
Hence
$B_{{\rm dist}(0,[p_{2},f_{11}(z_{0})])}(0)\subset S''$.
Combining this inclusion with (iii)' and (iv)', we have
\begin{equation}\label{Disc0}
f_{(1)^{n}}(z_{0})\in
B_{|f_{(1)^{N+2}}(z_{0})|}(0)\subset
B_{{\rm dist}(0,[p_{2},f_{11}(z_{0})])}(0)\subset S''
\end{equation}
for all $n$ with $n\ge N+2$.
By \eqref{Disc0}, if $n\ge N+3$, then all vertices $0$,
$f_{(1)^{n-1}}(z_{0})$, and $f_{(1)^{n}}(z_{0})$
of $\widetilde{A_{n}}$ belong to $S''$. Hence
\begin{equation}\label{inclusionH}
\widetilde{A_{n}}\subset S''\quad \mbox{for $n\ge N+3$.}
\end{equation}

Finally, we consider $\widetilde{A_{N+2}}$.
The vertices of $\widetilde{A_{N+2}}$ are $0$,
$f_{(1)^{N+1}}(z_{0})$, and $f_{(1)^{N+2}}(z_{0})$.
(i)' and Lemma \ref{lemma1} (v) give $f_{(1)^{N+1}}(z_{0})\in S''$.
\eqref{Disc0} gives $f_{(1)^{N+2}}(z_{0})\in S''$. Hence
\begin{equation}\label{inclusionM}
\widetilde{A_{N+2}}\subset S''.
\end{equation}

By \eqref{inclusionL}, \eqref{inclusionH}, and \eqref{inclusionM},
$\bigcup_{n=1}^{\infty}\widetilde{A_{n}}\subset S''$ is obtained.
\qed

\vspace{-5mm}

\begin{figure}[H]
\begin{center}
\caption{The set $f_{12}(S')$}
\label{f12Sdash}
\includegraphics[clip, width=7.5cm]{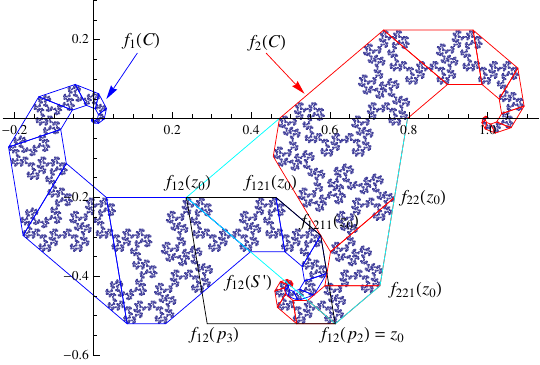}
\end{center}
\end{figure}

\vspace{-10mm}

\begin{figure}[H]
\begin{center}
\caption{The set $f_{22}(S')$}
\label{f22Sdash}
\includegraphics[clip, width=7.5cm]{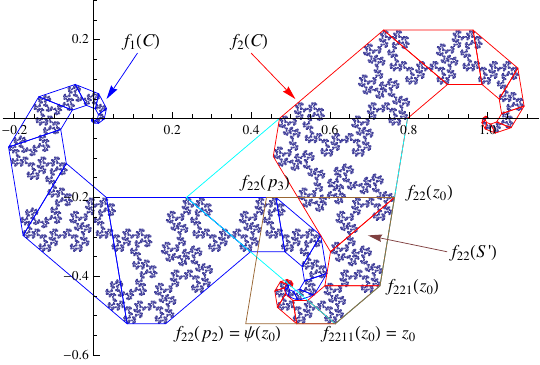}
\end{center}
\end{figure}

\vspace{-5mm}

\begin{proposition}\label{prop2}
Let $\xi$ be a fixed angle in $0<\xi<\pi/4$,
and $\mathcal{C}$ be the set of Definition \ref{def4}.
Then {\rm (i)} $f_{1}(\mathcal{C})\subset \mathcal{C}$ and
{\rm (ii)} $f_{2}(\mathcal{C})\subset \mathcal{C}$ are satisfied.
\end{proposition}

\noindent{\it Proof.}
(i)
It follows from the definitions of the sets $A_{n}$ and $\mathcal{C}$ that
\begin{equation*}
f_{1}(\mathcal{C})=\Big(\bigcup_{n=2}^{\infty}A_{n}\Big)
\cup f_{1}(B) \cup \Big(\bigcup_{n=1}^{\infty}f_{12}(A_{n})\Big).
\end{equation*}
Obviously, $\bigcup_{n=2}^{\infty}A_{n}\subset \mathcal{C}$.
Since $B\subset A_{0}$ by the definition of $B$,
$f_{1}(B)\subset A_{1}\subset \mathcal{C}$ holds.
Applying $f_{12}$ to the inclusion
$\bigcup_{n=1}^{\infty}A_{n}\subset S'$
of Lemma \ref{lemma3}, we have
\begin{equation*}
\bigcup_{n=1}^{\infty}f_{12}(A_{n})\subset f_{12}(S')
=\mathcal{P}(f_{12}(z_{0}),f_{121}(z_{0}),f_{1211}(z_{0}),
f_{12}(p_{2}),f_{12}(p_{3})),
\end{equation*}
see Figure \ref{f12Sdash}.
Hence $f_{1}(\mathcal{C})\subset \mathcal{C}$ is derived from
\begin{equation}\label{f12an}
\mathcal{P}(f_{12}(z_{0}),f_{121}(z_{0}),f_{1211}(z_{0}),f_{12}(p_{2}),f_{12}(p_{3}))
\subset \mathcal{C}.
\end{equation}

By \eqref{f12}, \eqref{z0}, and Lemma \ref{lemma1} (iii), we have
\begin{equation}
f_{12}(p_{2})=z_{0}, \label{f12p2}
\end{equation}
and, by \eqref{f12}, Lemma \ref{lemma1} (iv), and \eqref{f12p2} , we have
\begin{equation}
f_{12}(p_{3})=z_{0}-\frac{|\alpha|^{4}}{1-|\alpha|^{2}}. \label{f12p3}
\end{equation}
$f_{12}(p_{3})$ is \textcolor{red}{an internal division point of} $[z_{0},f_{1}(z_{0})]$.
In fact, by \eqref{f12p3} and \eqref{lower1},
\begin{align}
\angle (z_{0},f_{12}(p_{3}),f_{1}(z_{0}))
= \arg \Big(-1\times
\frac{1-|\alpha|^{2}-|\alpha|^{4}}{|\alpha|^{2}+|\alpha|^{4}}\Big)
= \arg (-1)=\pi.\label{lemma4deiru}
\end{align}
Define the set $\widetilde{B}$ by
$\widetilde{B}=\mathcal{P}
(z_{0},f_{12}(p_{3}),f_{12}(z_{0}),f_{2}(z_{0}),f_{212}(z_{0}),f_{221}(z_{0}))$.
$\widetilde{B}$ is a convex set, for which
$\angle (f_{221}(z_{0}),z_{0},f_{12}(p_{3}))=\theta+\xi$,
$\angle (z_{0},f_{12}(p_{3}),f_{12}(z_{0}))=\theta$,
$\angle (f_{12}(p_{3}),f_{12}(z_{0}),f_{2}(z_{0}))=3\xi$,
and the other angles are the same as those of $B$.
$\widetilde{B}$ is a subset of $\mathcal{C}$.
In fact, since $f_{12}(p_{3})$ belongs to $[z_{0},f_{1}(z_{0})]$
and $A_{1}$ is a convex set,
$\mathcal{P}(z_{0},f_{12}(p_{3}),f_{12}(z_{0}))$ is a subset of $A_{1}$,
hence $\widetilde{B}\subset A_{1}\cup B\subset \mathcal{C}$.
Thus $\widetilde{B}$ is a convex set contained in $\mathcal{C}$.

Now the inclusion \eqref{f12an} follows from
$f_{12}(p_{2})\in\widetilde{B}$,
$f_{121}(z_{0})\in\widetilde{B}$, and $f_{1211}(z_{0})\in\widetilde{B}$.
By \eqref{f12p2}, $f_{12}(p_{2})=z_{0}\in\widetilde{B}$.
Let us consider $f_{121}(z_{0})$.
Since both $f_{212}(z_{0})$ and $f_{221}(z_{0})$
belong to the boundary of $\widetilde{B}$, and,
by Remark \ref{rem1}, \eqref{f21}, \eqref{f12}, \eqref{z0}, and \eqref{lower1},
\begin{align}\label{f22naibun}
\angle (f_{212}(z_{0}),f_{22}(z_{0}),f_{221}(z_{0}))
= \arg \Big(-1\times
\frac{|\alpha|^{2}}{1-|\alpha|^{2}-|\alpha|^{4}}\Big)
= \arg (-1)=\pi
\end{align}
holds, we have $f_{22}(z_{0})\in[f_{212}(z_{0}),f_{221}(z_{0})])\subset\widetilde{B}$.
Since both $f_{12}(z_{0})$ and $f_{22}(z_{0})$
belong to the boundary of $\widetilde{B}$, and,
by \eqref{f12}, \eqref{f22}, and \eqref{z0},
\begin{align*}
\angle (f_{12}(z_{0}),f_{121}(z_{0}),f_{22}(z_{0}))
= \arg \Big(-1\times
\frac{1-|\alpha|^{2}}{|\alpha|^{2}}\Big)
= \arg (-1)=\pi
\end{align*}
holds, we have $f_{121}(z_{0})\in[f_{12}(z_{0}),f_{22}(z_{0})])
\subset\widetilde{B}$.
Let us consider $f_{1211}(z_{0})$.
Since both $f_{121}(z_{0})$ and $f_{221}(z_{0})$ belong to $\widetilde{B}$,
and, by \eqref{f12}, \eqref{f22}, and \eqref{z0} that
$$
\angle (f_{121}(z_{0}),f_{1211}(z_{0}),f_{221}(z_{0}))
= \arg \Big(-1\times
\frac{1-|\alpha|^{2}}{|\alpha|^{2}}\Big)
= \arg (-1)= \pi
$$
holds, we have
$f_{1211}(z_{0})\in [f_{121}(z_{0}),f_{221}(z_{0})]\subset\widetilde{B}$.
This completes the proof of (i).

(ii)
It follows from the definition of the sets $A_{n}$ and $\mathcal{C}$ that
\begin{equation*}
f_{2}(\mathcal{C})=\Big(\bigcup_{n=1}^{\infty}f_{2}(A_{n})\Big)
\cup f_{2}(B) \cup \Big(\bigcup_{n=1}^{\infty}f_{22}(A_{n})\Big).
\end{equation*}
Obviously, $\bigcup_{n=1}^{\infty}f_{2}(A_{n})\subset \mathcal{C}$.
From $B\subset A_{0}$ and \eqref{f2nobunnkai},
it follows that $f_{2}(B)\subset \psi(A_{1})
=\mathcal{P}(\psi(z_{0}),f_{2}(z_{0}),f_{212}(z_{0}),f_{22}(z_{0}))$.
Since both $z_{0}$ and $f_{12}(p_{3})$
belong to the boundary of $\widetilde{B}$, and,
by \eqref{defpsi}, \eqref{z0}, \eqref{f12p3},
\begin{align*}
\angle (z_{0},\psi(z_{0}),f_{12}(p_{3}))
= \arg (-1\times |\alpha|^{2})= \arg (-1)=\pi
\end{align*}
holds, we have
$\psi(z_{0})\in [z_{0},f_{12}(p_{3})]\subset\widetilde{B}$.
Therefore $f_{2}(B)\subset\psi(A_{1})\subset \widetilde{B}\subset \mathcal{C}$.
Applying $f_{22}$ to the inclusion
$\bigcup_{n=1}^{\infty}A_{n}\subset S'$
of Lemma \ref{lemma3}, we have
\begin{equation*}
\bigcup_{n=1}^{\infty}f_{22}(A_{n})\subset f_{22}(S')
=\mathcal{P}(f_{22}(z_{0}),f_{221}(z_{0}),f_{2211}(z_{0}),
f_{22}(p_{2}),f_{22}(p_{3})),
\end{equation*}
see Figure \ref{f22Sdash}.
Hence $f_{2}(\mathcal{C})\subset \mathcal{C}$ is derived from
\begin{equation}\label{f22an}
\mathcal{P}(f_{22}(z_{0}),f_{221}(z_{0}),f_{2211}(z_{0}),f_{22}(p_{2}),f_{22}(p_{3}))
\subset \mathcal{C}.
\end{equation}

Obviously, $f_{2211}(z_{0})=z_{0}\in\widetilde{B}$.
We have already seen $f_{22}(z_{0})\in\widetilde{B}$.
Applying $\psi$ to \eqref{f12p2}, we have $f_{22}(p_{2})=\psi(z_{0})$.
We have already seen $\psi (z_{0})\in\widetilde{B}$,
hence $f_{22}(p_{2})\in \widetilde{B}$.
Since both $f_{12}(z_{0})$ and $f_{22}(z_{0})$
belong to the boundary of $\widetilde{B}$, and,
by \eqref{f12}, \eqref{f22}, \eqref{z0}, \eqref{lower1}, and Lemma \ref{lemma1} (iv),
\begin{align*}
\angle (f_{12}(z_{0}),f_{22}(p_{3}),f_{22}(z_{0}))
= \arg \Big(-1\times
\frac{|\alpha|^{2}+|\alpha|^{4}}{1-|\alpha|^{2}-|\alpha|^{4}}\Big)
= \arg (-1)=\pi
\end{align*}
holds, we have
$f_{22}(p_{3})\in [f_{12}(z_{0}),f_{22}(z_{0})]\subset\widetilde{B}$.
Thus \eqref{f22an} is valid.
This completes the proof of (ii). \qed

\section{A sufficient condition for
$f_{1}(\mathcal{C})\cap f_{2}(\mathcal{C})=\emptyset$}\label{sufficientII}

To prove $f_{1}(\mathcal{C})\cap f_{2}(\mathcal{C})=\emptyset$,
we first discuss that
$f_{1}(\mathcal{C})$ restricted \textcolor{red}{around} the point $0$
and $f_{2}(\mathcal{C})$ restricted \textcolor{red}{around} the point $1$ are disjoint,
see Figure \ref{majiwari} below.
Then $f_{1}(\mathcal{C})\cap f_{2}(\mathcal{C})=\emptyset$
follows from the disjointness between $f_{1}(\mathcal{C})$
and $f_{2}(\mathcal{C})$ \textcolor{red}{around} the point $\alpha$,
see Proposition \ref{prop3} below.

\vspace{-5mm}

\begin{figure}[H]
\begin{center}
\caption{The set
$\mathcal{P}(f_{2}(z_{0}),f_{21}(z_{0}),f_{211}(z_{0}),f_{1}^{-1}(z_{0}))$}
\label{cover1}
\includegraphics[clip, width=7.cm]{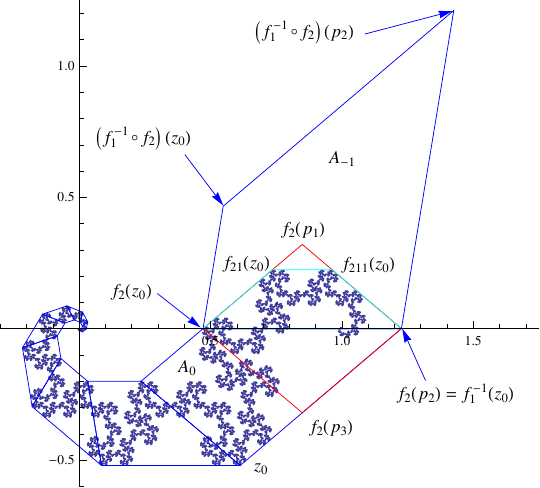}
\end{center}
\end{figure}

\begin{lemma}\label{lemma4}
For $ A_{-1}$ in Definition \ref{def3},
we have
$
\mathcal{P}(f_{2}(z_{0}),f_{21}(z_{0}),f_{211}(z_{0}),f_{1}^{-1}(z_{0}))
\subset A_{-1}
$, see Figure \ref{cover1}.
\end{lemma}

\noindent{\it Proof.}
Let $p_{2}$ be the point of Definition \ref{def5}.
Then, by \eqref{f12p2},
$A_{-1}$ of Definition \ref{def3} is rewritten as
$A_{-1}=\mathcal{P}(\textcolor{red}{(} f_{1}^{-1}\circ f_{2}\textcolor{red}{)} (p_{2}),
f_{2}(p_{2}),f_{2}(z_{0}),\textcolor{red}{(} f_{1}^{-1}\circ f_{2}\textcolor{red}{)} (z_{0}))$.
The sets $A_{1}$ and $A_{-1}$ are similar, hence, by Remark \ref{rem5},
$f_{2}(z_{0})$ is \textcolor{red}{an internal division point of} $[0,f_{2}(p_{2})]$.
Moreover, the half-line $HL(f_{2}(z_{0}),f_{2}(p_{1}))$
and the segment $[0,\textcolor{red}{(} f_{1}^{-1}\circ f_{2}\textcolor{red}{)} (p_{2})]$ are parallel, because
$\angle (f_{2}(p_{2}),f_{2}(z_{0}),f_{2}(p_{1}))
=\angle (p_{2},z_{0},p_{1})=\xi$ and
$\angle (f_{2}(p_{2}),0,\textcolor{red}{(} f_{1}^{-1}\circ f_{2}\textcolor{red}{)} (p_{2}))
=\angle (z_{0},0,f_{1}^{-1}(z_{0}))=\xi$.
Hence, if $f_{2}(p_{1})$ belongs to
$V^{+}(f_{2}(p_{2}),\textcolor{red}{(} f_{1}^{-1}\circ f_{2}\textcolor{red}{)} (p_{2}))$,
then $f_{2}(p_{1})$ belongs to $A_{-1}$.
In fact,
\begin{align*}
& \angle (\textcolor{red}{(} f_{1}^{-1}\circ f_{2}\textcolor{red}{)} (p_{2}),f_{2}(p_{2}),f_{2}(p_{1}))\\
& =\angle (\textcolor{red}{(} f_{1}^{-1}\circ f_{2}\textcolor{red}{)} (p_{2}),f_{2}(p_{2}),f_{2}(z_{0}))
-\angle (f_{2}(p_{1}),f_{2}(p_{2}),f_{2}(z_{0}))\\
& =\angle (f_{1}^{-1}(z_{0}),z_{0},f_{12}(z_{0}))
-\angle (p_{1},p_{2},z_{0})=\theta-\xi,
\end{align*}
where $\theta-\xi$ is greater than $0$ and less than $\pi$.
Thus $f_{2}(p_{1})\in A_{-1}$, and this gives
$[f_{2}(z_{0}),f_{2}(p_{1})]\subset A_{-1}$
and $[f_{2}(p_{1}),f_{2}(p_{2})]\subset A_{-1}$.
Applying $f_{2}$ to
$\angle (z_{0},f_{1}(z_{0}),p_{1})=\angle (p_{2},f_{11}(z_{0}),p_{1})=\pi$
in Remark \ref{rem6}, we have $f_{21}(z_{0}) \in [f_{2}(z_{0}),f_{2}(p_{1})]$
and $f_{211}(z_{0}) \in [f_{2}(p_{2}),f_{2}(p_{1})]$.
Hence, $f_{21}(z_{0}) \in A_{-1}$ and $f_{211}(z_{0}) \in A_{-1}$,
and all vertices of
$\mathcal{P}(f_{2}(z_{0}),f_{21}(z_{0}),f_{211}(z_{0}),f_{1}^{-1}(z_{0}))$
belong to $A_{-1}$.
\qed

\vspace{-3mm}

\begin{figure}[H]
\begin{center}
\caption{The sets $T$ and $T'$}
\label{cover2}
\includegraphics[clip, width=9.cm]{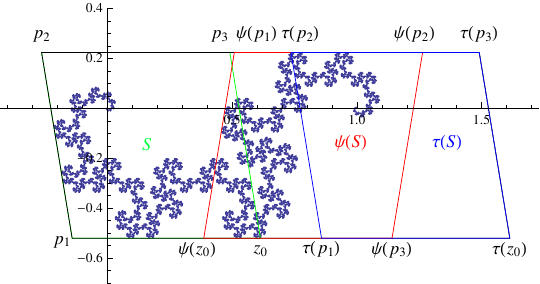}
\end{center}
\end{figure}

\begin{lemma}\label{lemma5}
Let $\xi$ be a fixed angle in $0<\xi<\pi/4$,
$\mathcal{C}$ be the set of Definition \ref{def4},
and $S$ be the set of Definition \ref{def5}.
Let $\psi$ be the function of \eqref{defpsi},
and define the function $\tau$ by $\tau(z)=z+1$.
Define the sets $T$ and $T'$ by
$T=\mathcal{P}(p_{1},p_{2},\tau(p_{3}),\tau(z_{0}))$
and $T'=\mathcal{P}(p_{1},p_{2},\psi(p_{2}),\psi(p_{3}))$,
see Figure \ref{cover2}. Then we have
{\rm (i)} $S\cup\psi (S)\subset T'$,
{\rm (ii)} $S\cup\tau (S)\subset T$,
{\rm (iii)} $T' \subset T$,
and {\rm (iv)} $\mathcal{C}\subset T'$.
\end{lemma}

\noindent{\it Proof.}
By \eqref{defpsi}, \eqref{z0}, and Lemma \ref{lemma1}, it is verified that
\begin{align}
& \psi(z_{0})=z_{0}-\frac{|\alpha|^{4}}{1-|\alpha|^{4}}, \label{psiz0}\\
& \psi(p_{1})=p_{2}+\frac{1-2|\alpha|^{4}}{1-|\alpha|^{4}}, \label{psip1}\\
& \psi(p_{2})=p_{2}+\frac{1+|\alpha|^{2}-|\alpha|^{4}}{1-|\alpha|^{4}}, \label{psip2}\\
& \psi(p_{3})=z_{0}+\frac{|\alpha|^{2}}{1-|\alpha|^{4}}. \label{psip3}
\end{align}

(i) \textcolor{red}{We get} $S\cup\psi (S)\subset T'$ from the statements that
both $\psi(z_{0})$ and $z_{0}$ belong to $[p_{1},\psi(p_{3})]$ and
both $p_{3}$ and $\psi(p_{1})$ belong to $[p_{2},\psi(p_{2})]$, \textcolor{red}{where}
$\psi(z_{0})\in [p_{1},\psi(p_{3})]$ is obtained by
Lemma \ref{lemma1} (ii), \eqref{psiz0}, and \eqref{psip3},  $z_{0}\in [p_{1},\psi(p_{3})]$ is obtained by
Lemma \ref{lemma1} (ii) and \eqref{psip3},
$p_{3}\in [p_{2},\psi(p_{2})]$ is obtained by
Lemma \ref{lemma1} (iv) and \eqref{psip2},
and $\psi(p_{1})\in [p_{2},\psi(p_{2})]$ is obtained by
\eqref{psip1} and \eqref{psip2}.

(ii) \textcolor{red}{We get} $S\cup\tau (S)\subset T$ from the statements that
both $z_{0}$ and $\tau(p_{1})$ belong to $[p_{1},\tau(z_{0})]$ and
both $p_{3}$ and $\tau(p_{2})$ belong to $[p_{2},\tau(p_{3})]$, \textcolor{red}{where}
$z_{0}, \tau(p_{1}) \in [p_{1},\tau(z_{0})]$
are obtained by Lemma \ref{lemma1} (ii) and the definition of $\tau$\textcolor{red}{, and}
$p_{3}, \tau(p_{2})\in [p_{2},\tau(p_{3})]$
are obtained by Lemma \ref{lemma1} (iv) and the definition of $\tau$.

(iii) \textcolor{red}{We get} $T'\subset T$ from the statements that
$\psi(p_{3})$ belongs to $[p_{1},\tau(z_{0})]$ and
$\psi(p_{2})$ belongs to $[p_{2},\tau(p_{3})]$, \textcolor{red}{where}
$\psi(p_{3})\in [p_{1},\tau(z_{0})]$ is obtained by
Lemma \ref{lemma1} (ii), \eqref{psip3}, and the definition of $\tau$\textcolor{red}{, and}
$\psi(p_{2})\in [p_{2},\tau(p_{3})]$ is obtained by
Lemma \ref{lemma1} (iv), \eqref{psip2}, and the definition of $\tau$.

(iv)
{\color{red}
From the expression
$
\mathcal{C}=\big(\bigcup_{n=1}^{\infty}A_{n}\big)
\cup B \cup \big(\bigcup_{n=1}^{\infty}f_{2}(A_{n})\big)
$,
it suffices to prove $\bigcup_{n=1}^{\infty}A_{n}\subset T'$,
$\bigcup_{n=1}^{\infty}f_{2}(A_{n})\subset T'$, and $B\subset T'$.
The inclusions $\bigcup_{n=1}^{\infty}A_{n}\subset S$ and
$S\subset S\cup\psi(S)\subset T'$ give
$\bigcup_{n=1}^{\infty}A_{n}\subset T'$.
}
Applying $\psi$ to the inclusion $\bigcup_{n=1}^{\infty}A_{n}\subset S$,
we have
$\bigcup_{n=1}^{\infty}\psi(A_{n})\subset \psi(S)$.
Hence, by \eqref{f2nobunnkai},
\begin{equation}\label{f2spiral}
\bigcup_{n=1}^{\infty}f_{2}(A_{n})\subset
\bigcup_{m=0}^{\infty}f_{2}(A_{m})
=\bigcup_{n=1}^{\infty}\psi(A_{n})
\subset \psi(S)\subset T'.
\end{equation}
Since $B$ is a subset of $A_{0}$,
{\color{red}
to show the inclusion $B\subset T'$,}
it \textcolor{red}{suffices} to prove that all vertices of $A_{0}$ belong to $T'$.
{\color{red}
By \eqref{f12p2}, $A_{0}$ is expressed as
$A_{0}=\mathcal{P}(f_{2}(p_{2}),z_{0},f_{12}(z_{0}),f_{2}(z_{0}))$.
}We have $z_{0}\in S\subset T'$.
Applying $\psi$ to $f_{1}(z_{0})\in [z_{0},p_{1}]$ (see Remark \ref{rem6}),
then using \eqref{f2nobunnkai}, we have
$f_{2}(z_{0})\in [\psi(z_{0}),\psi(p_{1})]$.
Hence $f_{2}(z_{0})\in \psi(S)\subset T'$.
Since $f_{12}(z_{0})\in [0,z_{0}]$ (see Remark \ref{rem5}),
we have $f_{12}(z_{0})\in S\subset T'$.
\textcolor{red}{By} \eqref{f2nobunnkai} and Lemma \ref{lemma1} (iii) (iv),
$$
\angle (\psi(p_{2}),f_{2}(p_{2}),\psi(p_{3}))
= \angle (p_{2},f_{1}(p_{2}),p_{3})
= \arg \Big(-1\times
\frac{1}{|\alpha|^{2}}\Big)= \arg (-1)=\pi
$$
holds, \textcolor{red}{and hence} $f_{2}(p_{2})\in [\psi(p_{2}),\psi(p_{3})]\subset\psi(S)\subset T'$.
\qed

\vspace{-5mm}

\begin{figure}[H]
\begin{center}
\caption{The sets $f_{1}(T)$, $f_{2}(T)$, and $f_{1}(T)\cap f_{2}(T)$}
\label{majiwari}
\includegraphics[clip, width=9.cm]{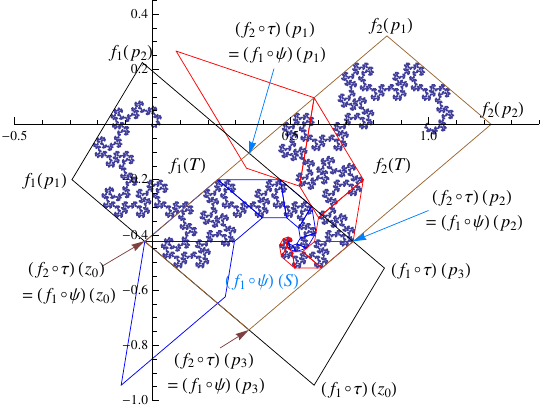}
\end{center}
\end{figure}

\begin{lemma}\label{lemma6}
Under the same \textcolor{red}{notation} as in Lemma \ref{lemma5}, we have
$f_{1}(T) \cap f_{2}(T)=\textcolor{red}{(} f_{1}\circ\psi\textcolor{red}{)} (S)$,
see Figure \ref{majiwari}.
\end{lemma}

\noindent{\it Proof.}
Applying $f_{2}$ to the expression
$T=(\mathcal{P}(p_{1},p_{2},\tau(p_{2}),\tau(p_{1}))
\setminus [\tau(p_{2}),\tau(p_{1})])
\cup\tau(S)$,
we have
$f_{2}(T)
=(\mathcal{P}(f_{2}(p_{1}),f_{2}(p_{2}),
\textcolor{red}{(} f_{2}\circ\tau\textcolor{red}{)} (p_{2}),
\textcolor{red}{(} f_{2}\circ\tau\textcolor{red}{)} (p_{1}))
\setminus [\textcolor{red}{(} f_{2}\circ\tau\textcolor{red}{)} (p_{2}),
\textcolor{red}{(} f_{2}\circ\tau\textcolor{red}{)} (p_{1})])
\cup \textcolor{red}{(} f_{2}\circ \tau\textcolor{red}{)} (S)$.
Then the statement is derived from

(i) $f_{1}(T) \cap
(\mathcal{P}(f_{2}(p_{1}),f_{2}(p_{2}),
\textcolor{red}{(} f_{2}\circ\tau\textcolor{red}{)} (p_{2}),
\textcolor{red}{(} f_{2}\circ\tau\textcolor{red}{)} (p_{1}))
\setminus [\textcolor{red}{(} f_{2}\circ\tau\textcolor{red}{)} (p_{2}),
\textcolor{red}{(} f_{2}\circ\tau\textcolor{red}{)} (p_{1})])
=\emptyset$,

(ii) $f_{1}(T) \cap \textcolor{red}{(} f_{2}\circ \tau\textcolor{red}{)} (S)
=\textcolor{red}{(}f_{1}\circ\psi\textcolor{red}{)} (S)$.

It is verified that
\begin{equation}
f_{2}\circ \tau=f_{1}\circ \psi. \label{kakan}
\end{equation}
By \eqref{kakan} and $\psi(S)\subset T$ (see Lemma \ref{lemma5} (i) (iii)),
$\textcolor{red}{(} f_{2}\circ \tau\textcolor{red}{)} (S)
=\textcolor{red}{(} f_{1}\circ\psi\textcolor{red}{)} (S)\subset f_{1}(T)$ holds,
which gives (ii).

Let us consider (i).
By \eqref{z0}, Lemma \ref{lemma1}, the definition of $\tau$,
and the inequality $1/2<|\alpha|<1/\sqrt{2}$, it is verified that
\begin{align*}
& \angle (\textcolor{red}{(} f_{1}\circ\tau\textcolor{red}{)} (p_{3}),
\textcolor{red}{(} f_{2}\circ\tau\textcolor{red}{)}  (p_{1}),f_{1}(p_{2}))
=\arg \Big(-1\times
\frac{1-2|\alpha|^{4}}{|\alpha|^{2}+2|\alpha|^{4}}\Big)
= \arg (-1)=\pi, \\
& \angle (\textcolor{red}{(} f_{1}\circ\tau\textcolor{red}{)} (p_{3}),
\textcolor{red}{(} f_{2}\circ\tau\textcolor{red}{)} (p_{2}),f_{1}(p_{2}))
=\arg \Big(-1\times
\frac{1+|\alpha|^{2}-|\alpha|^{4}}{|\alpha|^{4}}\Big)
= \arg (-1)=\pi, \\
& \angle (f_{2}(p_{1}),
\textcolor{red}{(}f_{2}\circ\tau\textcolor{red}{)} (p_{1}),f_{1}(p_{2}))
=\arg \Big(-\frac{\alpha}{\overline{\alpha}}\times
\frac{1+|\alpha|^{2}-|\alpha|^{4}}{1-|\alpha|^{4}}\Big)
= \arg \Big(-\frac{\alpha}{\overline{\alpha}}\Big)=\theta, \\
& \angle (f_{2}(p_{2}),\textcolor{red}{(} f_{2}\circ\tau\textcolor{red}{)} (p_{2}),f_{1}(p_{2}))
=\arg \Big(-\frac{\alpha}{\overline{\alpha}}\times
\frac{1-2|\alpha|^{4}}{1-|\alpha|^{4}}\Big)
= \arg \Big(-\frac{\alpha}{\overline{\alpha}}\Big)=\theta.
\end{align*}
Hence the points
$\textcolor{red}{(} f_{2}\circ\tau\textcolor{red}{)} (p_{1})$ and
$\textcolor{red}{(} f_{2}\circ\tau\textcolor{red}{)} (p_{2})$ are on
$L(f_{1}(p_{2}),\textcolor{red}{(} f_{1}\circ\tau\textcolor{red}{)} (p_{3}))$, and
the points $f_{2}(p_{1})$ and $f_{2}(p_{2})$ belong to
$V^{+}(f_{1}(p_{2}),\textcolor{red}{(} f_{1}\circ\tau\textcolor{red}{)}  (p_{3}))$.
Thus $\mathcal{P}(f_{2}(p_{1}),f_{2}(p_{2}),
\textcolor{red}{(} f_{2}\circ\tau\textcolor{red}{)} (p_{2}),
\textcolor{red}{(} f_{2}\circ\tau\textcolor{red}{)} (p_{1}))
\setminus [\textcolor{red}{(} f_{2}\circ\tau\textcolor{red}{)} (p_{2}),
\textcolor{red}{(} f_{2}\circ\tau\textcolor{red}{)} (p_{1})]$ is contained in
$V^{+}(f_{1}(p_{2}),\textcolor{red}{(} f_{1}\circ\tau\textcolor{red}{)} (p_{3}))$.
On the other hand, $f_{1}(T)$ is contained in
$\overline{V^{-}(f_{1}(p_{2}),\textcolor{red}{(} f_{1}\circ\tau\textcolor{red}{)} (p_{3}))}$.
Thus (i) is obtained. \qed

\vspace{-5mm}

\begin{figure}[H]
	\begin{center}
		\caption{The set $\mathcal{C}\cap \psi(S)$}
		\label{CpsiS}
		\includegraphics[clip, width=9.cm]{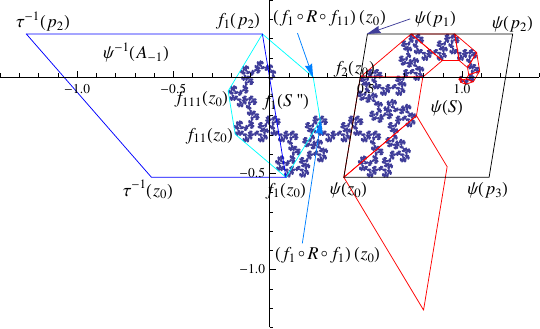}
	\end{center}
\end{figure}

\begin{lemma}\label{lemma7}
Let $\xi$ be a fixed angle in $0<\xi<\pi/4$,
$\mathcal{C}$ be the set of Definition \ref{def4},
$\psi$ be the function of \eqref{defpsi},
$S$ be the set of Definition \ref{def5},
and $A_{m}$ be the set of Definition \ref{def3}.
Then we have
$\mathcal{C}\cap \psi(S)
\subset \bigcup_{m=-1}^{\infty}f_{2}(A_{m})$,
see Figure \ref{CpsiS}.
\end{lemma}

\noindent{\it Proof.}
From \textcolor{red}{the expression}
$
\mathcal{C}=\big(\bigcup_{n=2}^{\infty}A_{n}\big)
\cup (A_{1}\cup B) \cup \big(\bigcup_{n=1}^{\infty}f_{2}(A_{n})\big)
$,
it \textcolor{red}{suffices} to prove the \textcolor{red}{following};

(i) $\bigcup_{n=1}^{\infty}f_{2}(A_{n})\subset \psi(S)$,

(ii) $(A_{1}\cup B) \cap \psi(S)\subset f_{2}(A_{0})\cup f_{2}(A_{-1})$,

(iii) $\Big(\bigcup_{n=2}^{\infty}A_{n}\Big) \cap \psi(S)=\emptyset$,

\noindent where (i) has been already proved as in \eqref{f2spiral}.

By \eqref{psiz0}, it is verified that
\begin{align*}
\angle (z_{0},\psi(z_{0}),f_{1}(z_{0}))
= \arg \Big(-1\times
\frac{1-|\alpha|^{2}}{|\alpha|^{2}}\Big)
= \arg (-1)=\pi,
\end{align*}
which means that $\psi(z_{0})$ belongs to the side $[z_{0},f_{1}(z_{0})]$ of $A_{1}$.
Hence $(A_{1}\cup B) \cap \psi(S)$ is expressed as
$\mathcal{P}(\psi(z_{0}),f_{2}(z_{0}),f_{212}(z_{0}),f_{221}(z_{0}),z_{0})$.
Since $f_{22}(z_{0})\in[f_{212}(z_{0}),f_{221}(z_{0})])$ (see \eqref{f22naibun}),
the above set is divided into two sets
$\mathcal{P}(\psi(z_{0}),f_{2}(z_{0}),f_{212}(z_{0}),f_{22}(z_{0}))$
and
$\mathcal{P}(f_{22}(z_{0}),f_{221}(z_{0}),z_{0},\psi(z_{0}))$,
where the former set is $f_{2}(A_{0})$ by \eqref{f2nobunnkai},
and the later set is contained in $f_{2}(A_{-1})$
by applying $f_{2}$ to the inclusion of Lemma \ref{lemma4}.
Thus (ii) is obtained.

Let us consider (iii).
By applying $f_{1}$ to the inclusion $\bigcup_{n=1}^{\infty}A_{n}\subset S''$,
where $S''$ is the set of Definition \ref{def5},
the set $\bigcup_{n=2}^{\infty}A_{n}$ is contained in $f_{1}(S'')$.
The set $f_{1}(S'')$ is divided into two sets
$\mathcal{P}(f_{1}(z_{0}),f_{11}(z_{0}),f_{111}(z_{0}),f_{1}(p_{2}))$
and
$\mathcal{P}(f_{1}(z_{0}),f_{1}(p_{2}),
\textcolor{red}{(} f_{1}\circ R \circ f_{11}\textcolor{red}{)} (z_{0}),
\textcolor{red}{(} f_{1}\circ R \circ f_{1}\textcolor{red}{)} (z_{0}))$,
where the former set is contained in $\psi^{-1}(A_{-1})$
by applying $\psi^{-1}$ to the inclusion of Lemma \ref{lemma4}.
Hence,
$$
\bigcup_{n=2}^{\infty}A_{n}\subset\psi^{-1}(A_{-1})\cup
\mathcal{P}(f_{1}(z_{0}),f_{1}(p_{2}),
\textcolor{red}{(} f_{1}\circ R \circ f_{11}\textcolor{red}{)} (z_{0}),
\textcolor{red}{(} f_{1}\circ R \circ f_{1}\textcolor{red}{)} (z_{0})),
$$
and it \textcolor{red}{suffices} for (iii) to prove the following (iv) and (v);

(iv) $\psi^{-1}(A_{-1}) \subset V^{+}(\psi(z_{0}),\psi(p_{1}))$,

(v) $\mathcal{P}(f_{1}(z_{0}),f_{1}(p_{2}),
\textcolor{red}{(}f_{1}\circ R \circ f_{11}\textcolor{red}{)} (z_{0}),
\textcolor{red}{(}f_{1}\circ R \circ f_{1}\textcolor{red}{)} (z_{0}))
\subset V^{+}(\psi(z_{0}),\psi(p_{1}))$.

By \eqref{f12p2} and \eqref{kakan},
the set $\psi^{-1}(A_{-1})$ in (iv) is expressed as
$$
\psi^{-1}(A_{-1})=
\mathcal{P}(\tau^{-1}(p_{2}),f_{1}(p_{2}), f_{1}(z_{0}),\tau^{-1}(z_{0})).
$$
By \eqref{defpsi}, \eqref{z0}, Lemma \ref{lemma1}, the definition of $\tau$,
and the inequality $1/2<|\alpha|<1/\sqrt{2}$,
it is verified that
\begin{align*}
& \angle (\psi(p_{1}),\psi(z_{0}),f_{1}(z_{0}))
=\arg \Big(-\frac{\alpha}{\overline{\alpha}}\times
\frac{1-|\alpha|^{2}}{1+|\alpha|^{2}}\Big)
= \arg \Big(-\frac{\alpha}{\overline{\alpha}}\Big)=\theta,\\
& \angle (f_{1}(p_{2}),\psi(p_{1}),\psi(z_{0}))
=\arg \Big(\frac{\overline{\alpha}}{\alpha}\times
\frac{|\alpha|^{2}+|\alpha|^{4}}{1-3|\alpha|^{4}}\Big)
= \arg \Big(\frac{\overline{\alpha}}{\alpha}\Big)=2\xi,\\
& \angle (\psi(p_{1}),\psi(z_{0}),\tau^{-1}(z_{0}))
=\arg \Big(-\frac{\alpha}{\overline{\alpha}}\times
\frac{1-|\alpha|^{2}}{|\alpha|^{2}}\Big)
= \arg \Big(-\frac{\alpha}{\overline{\alpha}}\Big)=\theta,\\
&
\angle (\tau^{-1}(p_{2}),\psi(p_{1}),\psi(z_{0}))
=\arg \Big(\frac{\overline{\alpha}}{\alpha}\times
\frac{|\alpha|^{2}+|\alpha|^{4}}{2-3|\alpha|^{4}}\Big)
= \arg \Big(\frac{\overline{\alpha}}{\alpha}\Big)=2\xi.
\end{align*}
Hence all vertices of $\psi^{-1}(A_{-1})$
belong to $V^{+}(\psi(z_{0}),\psi(p_{1}))$ and (iv) is obtained.

For (v), it \textcolor{red}{suffices} to prove that both
$\textcolor{red}{(} f_{1}\circ R \circ f_{11}\textcolor{red}{)} (z_{0})$ and
$\textcolor{red}{(} f_{1}\circ R \circ f_{1}\textcolor{red}{)} (z_{0})$
belong to $V^{+}(\psi(z_{0}),\psi(p_{1}))$.
The point $0$ belongs to $V^{+}(\psi(z_{0}),\psi(p_{1}))$,
because $\angle (f_{1}(z_{0}),0,f_{1}(p_{2}))=\pi$.
The point $f_{2}(z_{0})$ is on $L(\psi(z_{0}),\psi(p_{1}))$, because
\textcolor{red}{$$\angle (\psi(z_{0}),f_{2}(z_{0}),\psi(p_{1}))=\angle (z_{0},f_{1}(z_{0}),p_{1})=\pi$$}
(see Remark \ref{rem6}).
These give that the segment $[0,f_{2}(z_{0})]$
is contained in $\overline{V^{+}(\psi(z_{0}),\psi(p_{1}))}$.
By \eqref{z0}, $R(z)=(\alpha/\overline{\alpha})\overline{z}$,
and $1-2|\alpha|^{2}>0$, we have
\begin{align*}
\angle (0,\textcolor{red}{(} f_{1}\circ R \circ f_{11}\textcolor{red}{)} (z_{0}),f_{2}(z_{0}))
& = \arg \Big(-1\times
\frac{(1+|\alpha|^{2})(1-2|\alpha|^{2})}{|\alpha|^{4}}\Big)
= \arg (-1)=\pi.
\end{align*}
Thus, $\textcolor{red}{(} f_{1}\circ R \circ f_{11}\textcolor{red}{)} (z_{0})
\in [0,f_{2}(z_{0})]^{\circ}
\subset V^{+}(\psi(z_{0}),\psi(p_{1}))$.

Finally, we shall show that
$\textcolor{red}{(} f_{1}\circ R \circ f_{1}\textcolor{red}{)} (z_{0})$ belongs to $V^{+}(\psi(z_{0}),\psi(p_{1}))$.
It is verified that
$\textcolor{red}{(} f_{1}\circ R \circ f_{1}\textcolor{red}{)} (z_{0})
=|\alpha|^{2}z_{0}$, which means that
$\textcolor{red}{(} f_{1}\circ R \circ f_{1}\textcolor{red}{)} (z_{0})$
is on the segment $[0,z_{0}]$.
Let $\delta \in{\bf C}$ be the point such that
$\{ \delta \}=[0,z_{0}] \cap [f_{2}(z_{0}),\psi(z_{0})]$.
The triangles
$\mathcal{P}(\delta,0,f_{2}(z_{0}))$
and
$\mathcal{P}(\delta,z_{0},\psi(z_{0}))$
are similar. Since
\begin{equation*}
\Big|\frac{\psi(z_{0})-z_{0}}{f_{2}(z_{0})-0}\Big|
= \frac{|\alpha|^{4}}{1-|\alpha|^{2}-|\alpha|^{4}},
\end{equation*}
$|\delta|$ is expressed as
$|\delta|
= (1+\frac{|\alpha|^{4}}{1-|\alpha|^{2}-|\alpha|^{4}})^{-1}|z_{0}|$.
If the distance between $0$ and
$\textcolor{red}{(} f_{1}\circ R \circ f_{1}\textcolor{red}{)} (z_{0})$
is less than $|\delta|$, that is,
\begin{equation}\label{chikai}
|\alpha|^{2}|z_{0}|
<\frac{1}{1+\frac{|\alpha|^{4}}{1-|\alpha|^{2}-|\alpha|^{4}}}|z_{0}|,
\end{equation}
then $\textcolor{red}{(} f_{1}\circ R \circ f_{1}\textcolor{red}{)} (z_{0})$
belongs to the interior
of the segment $[0,\delta]$, and, consequently,
$\textcolor{red}{(} f_{1}\circ R \circ f_{1}\textcolor{red}{)} (z_{0})\in V^{+}(\psi(z_{0}),\psi(p_{1}))$.
The inequality \eqref{chikai} is equivalent to
the inequality $1-2|\alpha|^{2}>0$, which is valid
under the assumption $0<\xi<\pi/4$.
Thus $\textcolor{red}{(} f_{1}\circ R \circ f_{1}\textcolor{red}{)} (z_{0})
\in V^{+}(\psi(z_{0}),\psi(p_{1}))$.
\qed

\vspace{-6mm}

\begin{figure}[H]
	\begin{center}
		\caption{The set $\mathcal{C}\cap \tau(S)$}
		\label{CtauS}
		\includegraphics[clip, width=9.5cm]{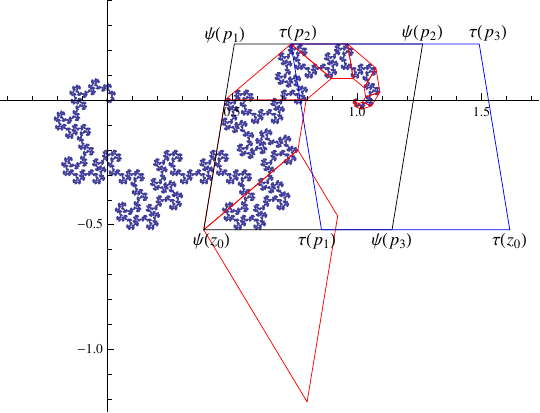}
	\end{center}
\end{figure}

\vspace{-5mm}

\begin{lemma}\label{lemma8}
Let $\xi$ be a fixed angle in $0<\xi<\pi/4$,
$\mathcal{C}$ be the set of Definition \ref{def4},
$\tau$ be the function of Lemma \ref{lemma5},
$S$ be the set of Definition \ref{def5},
and $A_{m}$ be the set of Definition \ref{def3}.
Then we have
$\mathcal{C}\cap \tau(S)
\subset \bigcup_{m=-1}^{\infty}f_{2}(A_{m})$,
see Figure \ref{CtauS}.
\end{lemma}

\noindent{\it Proof.}
Since $\mathcal{C}\subset T'$ (see Lemma \ref{lemma5} (iv)),
$\mathcal{C}\subset \overline{V^{+}(\psi(p_{3}),\psi(p_{2}))}$ holds.
Hence,
\begin{equation*}
\mathcal{C}\cap \tau(S)=
\mathcal{C}\cap (\tau(S) \cap \overline{V^{+}(\psi(p_{3}),\psi(p_{2}))}).
\end{equation*}
By \eqref{defpsi}, \eqref{z0}, Lemma \ref{lemma1}, and the definition of $\tau$,
it is verified that
\begin{align*}
\angle (\tau(p_{3}),\psi(p_{2}),\tau(p_{2}))
= \angle (\psi(p_{1}),\tau(p_{2}),\psi(p_{2}))
= \arg \Big(-1\times\frac{1}{|\alpha|^{2}}\Big)
= \arg (-1)=\pi.
\end{align*}
Hence $\psi(p_{2})$ (resp.\ $\tau(p_{2})$)
belongs to $[\tau(p_{2}),\tau(p_{3})]$ (resp.\ $[\psi(p_{1}),\psi(p_{2})]$).
If $[\tau(p_{2}),\tau(p_{1})]\cap[\psi(p_{2}),\psi(p_{3})]=\emptyset$,
then $\tau(p_{1})$ belongs to the interior of $[\psi(z_{0}),\psi(p_{3})]$, hence
\begin{equation*}
\tau(S) \cap \overline{V^{+}(\psi(p_{3}),\psi(p_{2}))}
= \mathcal{P}(\tau(p_{1}),\tau(p_{2}),\psi(p_{2}),\psi(p_{3}))
\subset \psi(S).
\end{equation*}
If $[\tau(p_{2}),\tau(p_{1})]\cap[\psi(p_{2}),\psi(p_{3})]\neq \emptyset$,
that is, $[\tau(p_{2}),\tau(p_{1})]\cap[\psi(p_{2}),\psi(p_{3})]=\{\varepsilon\}$, say,
then $\varepsilon$ belongs to $[\psi(p_{2}),\psi(p_{3}))]$, hence
\begin{equation*}
\tau(S) \cap \overline{V^{+}(\psi(p_{3}),\psi(p_{2}))}
= \mathcal{P}(\varepsilon,\tau(p_{2}),\psi(p_{2}))
\subset \psi(S).
\end{equation*}
Thus $\mathcal{C}\cap \tau(S) \subset \mathcal{C}\cap \psi(S)$ holds.
Combining this inclusion with Lemma \ref{lemma7}, we have this lemma.
\qed

\begin{proposition}\label{prop3}
Let $\xi$ be a fixed angle in $0<\xi<\pi/4$,
$\mathcal{C}$ be the set of Definition \ref{def4},
and $A_{m}$ be the set of Definition \ref{def3}.
If
\begin{equation}\label{majiwaranai}
\Big(\bigcup_{m=-1}^{\infty}f_{12}(A_{m})\Big) \cap
\Big(\bigcup_{m=-1}^{\infty}f_{22}(A_{m})\Big) =\emptyset
\end{equation}
is satisfied, then
$f_{1}(\mathcal{C})\cap f_{2}(\mathcal{C})=\emptyset$ is satisfied.
\end{proposition}

\noindent{\it Proof.}
By $\mathcal{C}\subset T$ (see Lemma \ref{lemma5}) and Lemma \ref{lemma6}, we have
$$
(f_{1}(\mathcal{C})\cap \textcolor{red}{(} f_{1}\circ\psi\textcolor{red}{)} (S))\cap
(f_{2}(\mathcal{C})\cap \textcolor{red}{(} f_{1}\circ\psi\textcolor{red}{)} (S))
= f_{1}(\mathcal{C})\cap f_{2}(\mathcal{C}).
$$
Applying $f_{1}$ to the inclusion of Lemma \ref{lemma7}, we have
$f_{1}(\mathcal{C})\cap \textcolor{red}{(} f_{1}\circ\psi\textcolor{red}{)} (S)
\subset \bigcup_{m=-1}^{\infty}f_{12}(A_{m})$.
Applying $f_{2}$ to the inclusion of Lemma \ref{lemma8},
then using \eqref{kakan}, we have
$f_{2}(\mathcal{C})\cap \textcolor{red}{(} f_{1}\circ\psi\textcolor{red}{)} (S)
\subset \bigcup_{m=-1}^{\infty}f_{22}(A_{m})$.
Hence, if \eqref{majiwaranai} is satisfied,
then $f_{1}(\mathcal{C})\cap f_{2}(\mathcal{C})=\emptyset$ is satisfied.
\qed

\section{Proof of $f_{1}(\mathcal{C})\cap f_{2}(\mathcal{C})=\emptyset$
for $0<\xi<\xi_{0}<\pi/4$}\label{proofII}

The aim of this section is to prove the following proposition.

\begin{proposition}\label{prop4}
Let $\mathcal{C}$ be the set of Definition \ref{def4}.
There exists a constant $\xi_{0}$, whose approximate value is
$\xi_{0}\approx 0.703858$,
such that $f_{1}(\mathcal{C})\cap f_{2}(\mathcal{C})=\emptyset$
is satisfied for any $\xi$ in $0<\xi<\xi_{0}$.
\end{proposition}

\begin{lemma}\label{lemma9}
Define the function $g$ by
$g(z)=\alpha(z-\alpha)+\alpha$.
{\color{red}
The function $g$ has two expressions $g=f_{12}\circ f_{1}\circ f_{12}^{-1}$
and $g=f_{22}\circ f_{1}\circ f_{22}^{-1}$.
For any integer $m$, define $g^{m}$ by
$g^{m}=f_{12}\circ f_{(1)^{m}}\circ f_{12}^{-1}=f_{22}\circ f_{(1)^{m}}\circ f_{22}^{-1}$,
where $f_{(1)^{m}}$ is the same one as in Definition \ref{def3}.
Then we have
{\rm (i)} $f_{12}\circ f_{(1)^{m}}=g^{m}\circ f_{12}$
and
{\rm (ii)} $f_{22}\circ f_{(1)^{m}}=g^{m}\circ f_{22}$.}
\end{lemma}

\noindent{\it Proof.}
{\color{red}
This lemma directly follows from the definition of $g^{m}$.
}
\qed
\vspace{3mm}

\begin{figure}
\begin{center}
\caption{The sets $\bigcup_{m=-2}^{\infty}\textcolor{red}{(} g^{m}\circ f_{12}\textcolor{red}{)} (A_{1})$
and $\bigcup_{m=-2}^{\infty}\textcolor{red}{(} g^{m}\circ f_{22}\textcolor{red}{)} (A_{1})$}
\label{disjoint1}
\includegraphics[clip, width=8.cm]{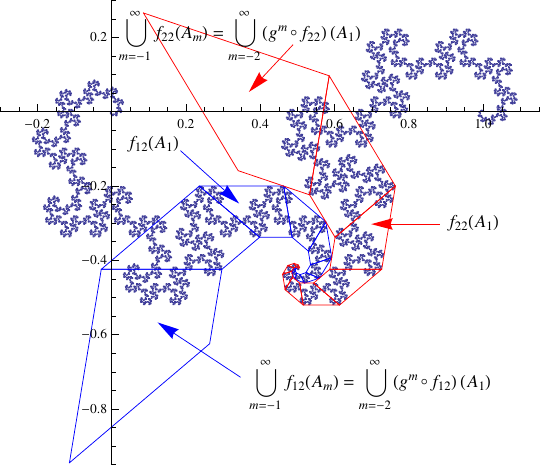}
\end{center}
\end{figure}

By Lemma \ref{lemma9}, $f_{12}(A_{m})$ and $f_{22}(A_{m})$ in the condition
\eqref{majiwaranai} are expressed as
$f_{12}(A_{m})=\textcolor{red}{(} g^{m-1}\circ f_{12}\textcolor{red}{)} (A_{1})$
and
$f_{22}(A_{m})=\textcolor{red}{(} g^{m-1}\circ f_{22}\textcolor{red}{)} (A_{1})$.
Hence, if
\begin{equation}\label{spiral00}
\Big(\bigcup_{m=-\infty}^{\infty}\textcolor{red}{(} g^{m}\circ f_{12}\textcolor{red}{)} (A_{1})\Big) \cap
\Big(\bigcup_{m=-\infty}^{\infty}\textcolor{red}{(} g^{m}\circ f_{22}\textcolor{red}{)} (A_{1})\Big) =\emptyset
\end{equation}
is satisfied, then \eqref{majiwaranai} is satisfied,
see Figure \ref{disjoint1}.

\begin{figure}
	\begin{center}
		\caption{${\rm Cone}(f_{12}(\widetilde{A_{1}}))$,
			${\rm Cone}(f_{22}(\widetilde{A_{1}}))$,
			and ${\rm Cone}(\textcolor{red}{(} \psi\circ f_{22}\textcolor{red}{)} (\widetilde{A_{1}}))$}
		\label{cones1}
		\includegraphics[clip, width=9.cm]{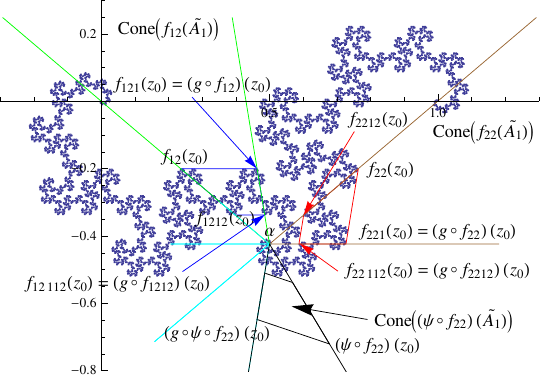}
	\end{center}
\end{figure}

Let $a$ be an element of
$\Big(\bigcup_{m=-\infty}^{\infty}\textcolor{red}{(} g^{m}\circ f_{12}\textcolor{red}{)} (A_{1})\Big) \cap
\Big(\bigcup_{m=-\infty}^{\infty}\textcolor{red}{(} g^{m}\circ f_{22}\textcolor{red}{)} (A_{1})\Big)$.
Then there exist integers $m$ and $l$ such that
$g^{-m}(a)\in f_{12}(A_{1})$ and $g^{-l}(a)\in f_{22}(A_{1})$.
In the case $m\ge l$, $g^{-l}(a)$ is replaced by $a'$ to give
$a'\in \textcolor{red}{(} g^{m-l}\circ f_{12}\textcolor{red}{)} (A_{1})\cap f_{22}(A_{1})$.
In the case $m\le l$, $g^{-m}(a)$ is replaced by $a'$ to give
$a'\in f_{12}(A_{1}) \cap \textcolor{red}{(} g^{l-m}\circ f_{22}\textcolor{red}{)} (A_{1})$.
Hence \eqref{spiral00} is equivalent to
\begin{equation}\label{spiral11}
\Big(\Big(\bigcup_{n\ge 0}\textcolor{red}{(} g^{n}\circ f_{12}\textcolor{red}{)} (A_{1})\Big)
\cap f_{22}(A_{1})\Big) \cup
\Big(f_{12}(A_{1})\cap
\Big(\bigcup_{n\ge 0}\textcolor{red}{(} g^{n}\circ f_{22}\textcolor{red}{)} (A_{1})\Big)\Big)=\emptyset.
\end{equation}
Define the cones ${\rm Cone}(f_{22}(\widetilde{A_{1}}))$
and ${\rm Cone}(f_{12}(\widetilde{A_{1}}))$ by
\begin{align*}
{\rm Cone}(f_{22}(\widetilde{A_{1}}))=
\overline{V^{-}(\alpha,f_{22}(z_{0}))} \cap
\overline{V^{+}(\alpha,\textcolor{red}{(} g\circ f_{22}\textcolor{red}{)} (z_{0}))},\\
{\rm Cone}(f_{12}(\widetilde{A_{1}}))=
\overline{V^{-}(\alpha,f_{12}(z_{0}))} \cap
\overline{V^{+}(\alpha,\textcolor{red}{(} g\circ f_{12}\textcolor{red}{)} (z_{0}))},
\end{align*}
see Figure \ref{cones1}. Then,
$f_{22}(A_{1})\subset f_{22}(\widetilde{A_{1}})
\subset {\rm Cone}(f_{22}(\widetilde{A_{1}}))$
and $f_{12}(A_{1})\subset f_{12}(\widetilde{A_{1}})
\subset {\rm Cone}(f_{12}(\widetilde{A_{1}}))$.
Hence the first (resp.\ second) set on the left-hand side of
\eqref{spiral11} is a subset of
$
\Big(\bigcup_{n\ge 0}\textcolor{red}{(} g^{n}\circ f_{12}\textcolor{red}{)} (\widetilde{A_{1}})\Big)
\cap {\rm Cone}(f_{22}(\widetilde{A_{1}}))\cap f_{22}(A_{1})
$
(resp.\
$
f_{12}(A_{1})\cap
\Big(\bigcup_{n\ge 0}\textcolor{red}{(} g^{n}\circ f_{22}\textcolor{red}{)} (\widetilde{A_{1}})\Big)
\cap {\rm Cone}(f_{12}(\widetilde{A_{1}}))
$
).
Thus \eqref{majiwaranai} is satisfied, if both
\begin{equation}\label{emptyA}
\Big(\bigcup_{n\ge 0}\textcolor{red}{(} g^{n}\circ f_{12}\textcolor{red}{)} (\widetilde{A_{1}})\Big)
\cap {\rm Cone}(f_{22}(\widetilde{A_{1}}))
\cap f_{22}(A_{1})=\emptyset
\end{equation}
and
\begin{equation}\label{emptyB}
f_{12}(A_{1})\cap
\Big(\bigcup_{n\ge 0}\textcolor{red}{(} g^{n}\circ f_{22}\textcolor{red}{)} (\widetilde{A_{1}})\Big)
\cap {\rm Cone}(f_{12}(\widetilde{A_{1}})) =\emptyset
\end{equation}
are satisfied.

Let us discuss \eqref{emptyA}.
Let $N(f_{12},f_{22})$ be the set of $k\ge 0$ such that
$\textcolor{red}{(} g^{k}\circ f_{12}\textcolor{red}{)} (z_{0})$ belongs to
${\rm Cone}(f_{22}(\widetilde{A_{1}}))\setminus
HL(\alpha,\textcolor{red}{(} g\circ f_{22}\textcolor{red}{)} (z_{0}))$.
For $k\in N(f_{12},f_{22})$,
the half-line $HL(\alpha,f_{22}(z_{0}))$
(resp.\ $HL(\alpha,\textcolor{red}{(} g\circ f_{22}\textcolor{red}{)} (z_{0}))$)
and the segment $[\textcolor{red}{(} g^{k-1}\circ f_{12}\textcolor{red}{)} (z_{0}),
\textcolor{red}{(} g^{k}\circ f_{12}\textcolor{red}{)} (z_{0})]$
(resp.\ $[\textcolor{red}{(} g^{k}\circ f_{12}\textcolor{red}{)} (z_{0}),
\textcolor{red}{(} g^{k+1}\circ f_{12}\textcolor{red}{)} (z_{0})]$)
are \textcolor{red}{intersecting}.
Define the points $\beta_{k}$ and $\gamma_{k}$ by
\begin{align*}
& \{\beta_{k}\}= HL(\alpha,f_{22}(z_{0}))
\cap [\textcolor{red}{(} g^{k-1}\circ f_{12}\textcolor{red}{)} (z_{0}),
\textcolor{red}{(} g^{k}\circ f_{12}\textcolor{red}{)} (z_{0})], \\
& \{\gamma_{k}\}= HL(\alpha,\textcolor{red}{(} g\circ f_{22}\textcolor{red}{)} (z_{0}))
\cap [\textcolor{red}{(} g^{k}\circ f_{12}\textcolor{red}{)} (z_{0}),
\textcolor{red}{(} g^{k+1}\circ f_{12}\textcolor{red}{)} (z_{0})].
\end{align*}
Then \eqref{emptyA} is equivalent to
\begin{equation}\label{emptyAdash}
	\Big(\bigcup_{k\in N(f_{12},f_{22})}
	\mathcal{P}(\alpha,\beta_{k}, \textcolor{red}{(} g^{k}\circ f_{12}\textcolor{red}{)} (z_{0}), \gamma_{k})\Big)
	\cap f_{22}(A_{1})=\emptyset,
\end{equation}
where
$\mathcal{P}(\alpha,\beta_{k}, \textcolor{red}{(} g^{k}\circ f_{12}\textcolor{red}{)} (z_{0}),\gamma_{k})$
is regarded as the triangle $\mathcal{P}(\alpha,\beta_{k}, \gamma_{k})$
in the case $\textcolor{red}{(} g^{k}\circ f_{12}\textcolor{red}{)} (z_{0})=\beta_{k}$.
By applying $g$ to the definition of $\beta_{k}$, and by $g(\alpha)=\alpha$,
$$
\{g(\beta_{k})\}=
HL(\alpha,\textcolor{red}{(} g\circ f_{22}\textcolor{red}{)} (z_{0}))
\cap [\textcolor{red}{(} g^{k}\circ f_{12}\textcolor{red}{)} (z_{0}),
\textcolor{red}{(} g^{k+1}\circ f_{12}\textcolor{red}{)} (z_{0})],
$$
hence $g(\beta_{k})=\gamma_{k}$. Then we have
\begin{align}
	\angle(\alpha,\beta_{k},\gamma_{k})
	= \angle(\alpha,\beta_{k},g(\beta_{k}))
	=\arg \overline{\alpha}=\xi. \label{kakudo1}
\end{align}
Hence $\mathcal{P}(\alpha,\beta_{k},\gamma_{k})$
is an isosceles triangle whose angles are
$\angle (\alpha,\beta_{k},\gamma_{k})
=\angle (\gamma_{k},\alpha,\beta_{k})=\xi$ and
$\angle (\beta_{k},\gamma_{k},\alpha)=\theta$.
Let $(k_{j})_{j=1}^{\infty}$ be the sequence such that
$k_{j}\in N(f_{12},f_{22})$ and $k_{1}<k_{2}<k_{3}<\cdots$.
By \eqref{kakudo1}, $[\beta_{k_{j}},\gamma_{k_{j}}]$,
$j=1,2,3,\ldots$, are parallel.
Since $2\pi -\xi >7\xi$,
the inequality $k_{j+1}-k_{j}\ge 7$ holds for any $j$.
Then, by $g(\alpha)=\alpha$, $k_{j+1}-k_{j}-2\ge 5$, and $1/|\alpha|>\sqrt{2}$,
\begin{align}
|\gamma_{k_{j}}-\alpha|	
& \ge |\textcolor{red}{(} g^{k_{j}+1}\circ f_{12}\textcolor{red}{)} (z_{0})-\alpha| \nonumber \\
& = |\textcolor{red}{(} g^{k_{j}+2-k_{j+1}}\circ g^{k_{j+1}-1}\circ f_{12}\textcolor{red}{)} (z_{0})
-g^{k_{j}+2-k_{j+1}}(\alpha)| \nonumber \\
& = (1/|\alpha|)^{k_{j+1}-k_{j}-2}\times
|\textcolor{red}{(} g^{k_{j+1}-1}\circ f_{12}\textcolor{red}{)} (z_{0})-\alpha| \nonumber \\
& > (\sqrt{2})^{5}\times
|\textcolor{red}{(} g^{k_{j+1}-1}\circ f_{12}\textcolor{red}{)} (z_{0})-\alpha|
> |\textcolor{red}{(} g^{k_{j+1}-1}\circ f_{12}\textcolor{red}{)} (z_{0})-\alpha| \nonumber \\
& \ge |\beta_{k_{j+1}}-\alpha| \ge
|\textcolor{red}{(} g^{k_{j+1}}\circ f_{12}\textcolor{red}{)} (z_{0})-\alpha|.
\label{return1} 	
\end{align}
The inequality
$|\gamma_{k_{j}}-\alpha| > |\textcolor{red}{(} g^{k_{j+1}}\circ f_{12}\textcolor{red}{)} (z_{0})-\alpha|$
in \eqref{return1} shows that $\textcolor{red}{(} g^{k_{j+1}}\circ f_{12}\textcolor{red}{)} (z_{0})$
belongs to the ball $B_{|\gamma_{k_{j}}-\alpha|}(\alpha)$.
Hence $\textcolor{red}{(} g^{k_{j+1}}\circ f_{12}\textcolor{red}{)} (z_{0})$ belongs to
${\rm Cone}(f_{22}(\widetilde{A_{1}})) \cap
B_{|\gamma_{k_{j}}-\alpha|}(\alpha)$.
This and
$\angle(\beta_{k_{j}},\gamma_{k_{j}},\alpha)=\theta > \pi/2$
give that $\textcolor{red}{(} g^{k_{j+1}}\circ f_{12}\textcolor{red}{)} (z_{0})$ belongs to the triangle
$\mathcal{P}(\alpha,\beta_{k_{j}},\gamma_{k_{j}})$.
By combining inequalities
$|\gamma_{k_{j}}-\alpha| > |\beta_{k_{j+1}}-\alpha|$ in \eqref{return1} with
$|\beta_{k_{j+1}}-\alpha| >|\gamma_{k_{j+1}}-\alpha|$,
the sequence $(|\gamma_{k_{j}}-\alpha|)_{j=1}^{\infty}$
is strictly monotonically decreasing as $j\to\infty$.
Then, by recalling all $[\beta_{k_{j}},\gamma_{k_{j}}]$s are parallel,
$\beta_{k_{j}}$ and $\gamma_{k_{j}}$, $j=1,2,3,\ldots$, belong to the triangle
$\mathcal{P}(\alpha,\beta_{k_{1}},\gamma_{k_{1}})$. Hence
$\bigcup_{k\in N(f_{12},f_{22})}
\mathcal{P}(\alpha,\beta_{k},
\textcolor{red}{(} g^{k}\circ f_{12}\textcolor{red}{)} (z_{0}), \gamma_{k})$
in \eqref{emptyAdash} is \textcolor{red}{identical to}
$\mathcal{P}(\alpha,\beta_{k_{1}},
\textcolor{red}{(} g^{k_{1}}\circ f_{12}\textcolor{red}{)} (z_{0}), \gamma_{k_{1}})$.
Thus \eqref{emptyAdash} is equivalent to
\begin{equation}\label{emptyAdash2}
\mathcal{P}(\alpha,\beta_{k_{1}}, \textcolor{red}{(} g^{k_{1}}\circ f_{12}\textcolor{red}{)} (z_{0}), \gamma_{k_{1}})	\cap f_{22}(A_{1})=\emptyset.
\end{equation}

By Definition \ref{def3} and Lemma \ref{lemma9} (ii),
$f_{22}(A_{1})$ is expressed as
$\mathcal{P}(f_{22}(z_{0}),\textcolor{red}{(} g\circ f_{22}\textcolor{red}{)} (z_{0}),
\textcolor{red}{(} g\circ f_{2212}\textcolor{red}{)} (z_{0}),f_{2212}(z_{0}))$.
The segments $[f_{2212}(z_{0}),\textcolor{red}{(} g\circ f_{2212}\textcolor{red}{)} (z_{0})]$ and
$[\beta_{k_{1}},\gamma_{k_{1}}]$ are parallel.
In fact, $\angle(\alpha,\beta_{k_{1}},\gamma_{k_{1}})=\xi$ as in \eqref{kakudo1},
and, by Lemma \ref{lemma9} (ii) and Remark \ref{rem5},
$
\angle(\alpha,f_{2212}(z_{0}),\textcolor{red}{(} g\circ f_{2212}\textcolor{red}{)} (z_{0}))
= \angle(0,f_{12}(z_{0}),f_{112}(z_{0}))=\xi.
$

We shall show that the point $\textcolor{red}{(} g^{k_{1}}\circ f_{12}\textcolor{red}{)} (z_{0})$
belongs to $\overline{V^{+}(\beta_{k_{1}},\gamma_{k_{1}})}$.
If $\textcolor{red}{(} g^{k_{1}}\circ f_{12}\textcolor{red}{)} (z_{0})$
is \textcolor{red}{identical to} $\beta_{k_{1}}$,
then the statement is obvious.
Assume that $\textcolor{red}{(} g^{k_{1}}\circ f_{12}\textcolor{red}{)} (z_{0})$
is not \textcolor{red}{identical to} $\beta_{k_{1}}$.
Then the angle
$\angle (\textcolor{red}{(} g^{k_{1}}\circ f_{12}\textcolor{red}{)} (z_{0}),\alpha,\beta_{k_{1}})$ satisfies
$0<\angle (\textcolor{red}{(} g^{k_{1}}\circ f_{12}\textcolor{red}{)} (z_{0}),\alpha,\beta_{k_{1}})<\xi$.
This and
$\angle (\beta_{k_{1}},\textcolor{red}{(} g^{k_{1}}\circ f_{12}\textcolor{red}{)} (z_{0}),\alpha)=\theta$
give that $\xi<
\angle (\alpha,\beta_{k_{1}},\textcolor{red}{(} g^{k_{1}}\circ f_{12}\textcolor{red}{)} (z_{0}))<2\xi$,
hence $0<\angle(\alpha,\beta_{k_{1}},\gamma_{k_{1}})
<\angle (\alpha,\beta_{k_{1}},\textcolor{red}{(} g^{k_{1}}\circ f_{12}\textcolor{red}{)} (z_{0}))<\pi/2$.
Thus $\angle (\gamma_{k_{1}},\beta_{k_{1}}, \textcolor{red}{(} g^{k_{1}}\circ f_{12}\textcolor{red}{)} (z_{0}))$
is greater than $0$ and less than $\pi/2$, that is,
$\textcolor{red}{(} g^{k_{1}}\circ f_{12}\textcolor{red}{)} (z_{0})\in V^{+}(\beta_{k_{1}},\gamma_{k_{1}})$.

Next, we shall show that if the point $\textcolor{red}{(} g^{k_{1}}\circ f_{12}\textcolor{red}{)} (z_{0})$
belongs to $V^{+}(\textcolor{red}{(} g\circ f_{2212}\textcolor{red}{)} (z_{0}),f_{2212}(z_{0}))$,
then both $\beta_{k_{1}}$ and $\gamma_{k_{1}}$ \textcolor{red}{belong to}
$V^{+}(\textcolor{red}{(} g\circ f_{2212}\textcolor{red}{)} (z_{0}),f_{2212}(z_{0}))$.
If $\textcolor{red}{(} g^{k_{1}}\circ f_{12}\textcolor{red}{)} (z_{0})$
is \textcolor{red}{identical to} $\beta_{k_{1}}$,
then the statement is obvious.
Assume that $\textcolor{red}{(} g^{k_{1}}\circ f_{12}\textcolor{red}{)} (z_{0})$
is not identical \textcolor{red}{to} $\beta_{k_{1}}$.
Then we have already seen that $\textcolor{red}{(} g^{k_{1}}\circ f_{12}\textcolor{red}{)} (z_{0})$
belongs to $V^{+}(\beta_{k_{1}},\gamma_{k_{1}})$.
Hence, if $\textcolor{red}{(} g^{k_{1}}\circ f_{12}\textcolor{red}{)} (z_{0})$
belongs to $V^{+}(\textcolor{red}{(} g\circ f_{2212}\textcolor{red}{)} (z_{0}),f_{2212}(z_{0}))$,
then $\textcolor{red}{(} g^{k_{1}}\circ f_{12}\textcolor{red}{)} (z_{0})\in V^{+}(\beta_{k_{1}},\gamma_{k_{1}})
\cap V^{+}(\textcolor{red}{(} g\circ f_{2212}\textcolor{red}{)} (z_{0}),f_{2212}(z_{0}))$.
Let $\beta_{k_{1}}'$ and $\gamma_{k_{1}}'$ be the points such that
\begin{align*}
& \{\beta_{k_{1}}'\}=
L(\beta_{k_{1}},\textcolor{red}{(} g^{k_{1}}\circ f_{12}\textcolor{red}{)} (z_{0}))
\cap L(\textcolor{red}{(} g\circ f_{2212}\textcolor{red}{)} (z_{0}),f_{2212}(z_{0})), \\
& \{\gamma_{k_{1}}'\}=
L(\gamma_{k_{1}},\textcolor{red}{(} g^{k_{1}}\circ f_{12}\textcolor{red}{)} (z_{0}))
\cap L(\textcolor{red}{(} g\circ f_{2212}\textcolor{red}{)} (z_{0}),f_{2212}(z_{0})).
\end{align*}
Since the segments $[f_{2212}(z_{0}),\textcolor{red}{(} g\circ f_{2212}\textcolor{red}{)} (z_{0})]$ and
$[\beta_{k_{1}},\gamma_{k_{1}}]$ are parallel,
the triangles
$\mathcal{P}(\beta_{k_{1}}, \textcolor{red}{(} g^{k_{1}}\circ f_{12}\textcolor{red}{)} (z_{0}), \gamma_{k_{1}})$
and
$\mathcal{P}(\beta_{k_{1}}', \textcolor{red}{(} g^{k_{1}}\circ f_{12}\textcolor{red}{)} (z_{0}), \gamma_{k_{1}}')$
are similar. Hence $\beta_{k_{1}}$ (resp. $\gamma_{k_{1}}$)
is an \textcolor{red}{external division point of} $[\beta_{k_{1}}',\textcolor{red}{(} g^{k_{1}}\circ f_{12}\textcolor{red}{)} (z_{0})]$
(resp. $[\gamma_{k_{1}}',\textcolor{red}{(} g^{k_{1}}\circ f_{12}\textcolor{red}{)} (z_{0})]$)
Thus, if  $\textcolor{red}{(} g^{k_{1}}\circ f_{12}\textcolor{red}{)} (z_{0})$ belongs to
$V^{+}(\textcolor{red}{(} g\circ f_{2212}\textcolor{red}{)} (z_{0}),f_{2212}(z_{0}))$,
then both $\beta_{k_{1}}$ and $\gamma_{k_{1}}$
belong to $V^{+}(\textcolor{red}{(} g\circ f_{2212}\textcolor{red}{)} (z_{0}),f_{2212}(z_{0}))$.
The point $\alpha$ belongs to $V^{+}(\textcolor{red}{(} g\circ f_{2212}\textcolor{red}{)} (z_{0}),f_{2212}(z_{0}))$,
because
$\angle(f_{2212}(z_{0}),\textcolor{red}{(} g\circ f_{2212}\textcolor{red}{)} (z_{0}),\alpha)
= \angle(f_{12}(z_{0}),f_{112}(z_{0}),0)=\theta$.
Thus, if $\textcolor{red}{(} g^{k_{1}}\circ f_{12}\textcolor{red}{)} (z_{0})$ belongs to
$V^{+}(\textcolor{red}{(} g\circ f_{2212}\textcolor{red}{)} (z_{0}),f_{2212}(z_{0}))$,
$\mathcal{P}(\alpha,\beta_{k_{1}}, \textcolor{red}{(} g^{k_{1}}\circ f_{12}\textcolor{red}{)} (z_{0}), \gamma_{k_{1}})$
is contained in $V^{+}(\textcolor{red}{(} g\circ f_{2212}\textcolor{red}{)} (z_{0}),f_{2212}(z_{0}))$.
On the other hand, $f_{22}(A_{1})$ is contained in
$\overline{V^{-}(\textcolor{red}{(} g\circ f_{2212}\textcolor{red}{)} (z_{0}),f_{2212}(z_{0}))}$.
Hence, \eqref{emptyAdash2} (consequently \eqref{emptyA}) is satisfied if
$\textcolor{red}{(} g^{k_{1}}\circ f_{12}\textcolor{red}{)} (z_{0})
\in V^{+}(\textcolor{red}{(} g\circ f_{2212}\textcolor{red}{)} (z_{0}),f_{2212}(z_{0}))$
can be shown.

\textcolor{red}{A similar discussion} works on \eqref{emptyB},
and \textcolor{red}{details} in this case are omitted.
We have obtained the following.

\begin{lemma}\label{lemma10}
Let $N(f_{12},f_{22})$ \rm{(resp.\ }$N(f_{22},f_{12})$\rm{)}
be the set of $k\ge 0$ \rm{(resp.\ }$l\ge 0$\rm{)} such that
$\textcolor{red}{(} g^{k}\circ f_{12}\textcolor{red}{)} (z_{0})$
\rm{(resp.\ }$\textcolor{red}{(} g^{l}\circ f_{22}\textcolor{red}{)} (z_{0})$\rm{)} belongs to
${\rm Cone}(f_{22}(\widetilde{A_{1}}))\setminus HL(\alpha,\textcolor{red}{(} g\circ f_{22}\textcolor{red}{)} (z_{0}))$
\rm{(resp.\ ${\rm Cone}(f_{12}(\widetilde{A_{1}}))\setminus
HL(\alpha,\textcolor{red}{(} g\circ f_{12}\textcolor{red}{)} (z_{0}))$}\rm{)}.
Let $k_{1}$ \rm{(resp.\ }$l_{1}$\rm{)}
be the least number of $N(f_{12},f_{22})$
\rm{(resp.\ }$N(f_{22},f_{12})$\rm{)}.
If both \textcolor{red}{conditions}
{\rm (i)} $\textcolor{red}{(} g^{k_{1}}\circ f_{12}\textcolor{red}{)} (z_{0})
\in V^{+}(\textcolor{red}{(} g\circ f_{2212}\textcolor{red}{)} (z_{0}),f_{2212}(z_{0}))$
and
{\rm (ii)} $\textcolor{red}{(} g^{l_{1}}\circ f_{22}\textcolor{red}{)} (z_{0})
\in V^{+}(\textcolor{red}{(} g\circ f_{1212}\textcolor{red}{)} (z_{0}),f_{1212}(z_{0}))$
are satisfied, then
$f_{1}(\mathcal{C})\cap f_{2}(\mathcal{C})=\emptyset$ is satisfied.
\end{lemma}

\begin{lemma}\label{hidarigawa}
Let $\xi$ be a fixed angle in $0<\xi<\pi/4$.
Choose a positive integer $N\ge 3$ so that
the inequality $\pi/(N+2)\le \xi< \pi/(N+1)$ is satisfied.
If both \textcolor{red}{conditions}
{\rm (i)} $\textcolor{red}{(} g^{N}\circ f_{12}\textcolor{red}{)} (z_{0})
\in V^{+}(\textcolor{red}{(} g\circ f_{2212}\textcolor{red}{)} (z_{0}),f_{2212}(z_{0}))$
and
{\rm (ii)} $\textcolor{red}{(} g^{N+4}\circ f_{22}\textcolor{red}{)} (z_{0})
\in V^{+}(\textcolor{red}{(} g\circ f_{1212}\textcolor{red}{)} (z_{0}),f_{1212}(z_{0}))$
are satisfied, then
$f_{1}(\mathcal{C})\cap f_{2}(\mathcal{C})=\emptyset$ is satisfied.
\end{lemma}

\noindent{\it Proof.}
$f_{12}(z_{0})$  belongs to $V^{+}(\alpha,f_{22}(z_{0}))$,
because
$
\angle (f_{22}(z_{0}),\alpha,f_{12}(z_{0}))
= \arg (-\frac{\alpha}{\overline{\alpha}})=\theta
$
by \eqref{f12} and \eqref{f22}.
Applying $f_{12}$ to
$f_{(1)^{n}}(z_{0})\in W\setminus [0,q]$ (see Lemma \ref{lemma2} (i)),
then using Lemma \ref{lemma9} (i), we have
$\textcolor{red}{(} g^{n}\circ f_{12}\textcolor{red}{)} (z_{0})\in f_{12}(W\setminus [0,q])$
for all $n$ with $1\le n\le N-1$.
The equality $f_{12}(q)=f_{22}(z_{0})$ is verified by
\eqref{f12}, Lemma \ref{lemma1} (i), and \eqref{f22}. Hence,
$$
\textcolor{red}{(} g^{n}\circ f_{12}\textcolor{red}{)} (z_{0})\in
\mathcal{P}(\alpha,f_{12}(z_{0}),f_{22}(z_{0}))\setminus [\alpha,f_{22}(z_{0})]
\subset V^{+}(\alpha,f_{22}(z_{0}))
$$
for $1\le n\le N-1$.
Thus $\textcolor{red}{(} g^{n}\circ f_{12}\textcolor{red}{)} (z_{0})$ belongs to $V^{+}(\alpha,f_{22}(z_{0}))$
for all $n$ with $0\le n\le N-1$.

Let us show that $\textcolor{red}{(} g^{N}\circ f_{12}\textcolor{red}{)} (z_{0})$ belongs to
${\rm Cone}(f_{22}(\widetilde{A_{1}}))\setminus HL(\alpha,\textcolor{red}{(} g\circ f_{22}\textcolor{red}{)} (z_{0}))$.
Applying $f_{12}$ to $\angle (f_{(1)^{N}}(z_{0}),0,q)\ge 0$
and $\angle (-1,0,f_{(1)^{N}}(z_{0}))>0$ (see Lemma \ref{lemma2} (ii)),
then using Lemma \ref{lemma9} (i) and $f_{12}(q)=f_{22}(z_{0})$, we have
\begin{equation}\label{ConeUpper}
\angle (\textcolor{red}{(} g^{N}\circ f_{12}\textcolor{red}{)} (z_{0}),\alpha,f_{22}(z_{0}))\ge 0
\end{equation}
and
$
\angle (f_{12}(-1),\alpha, \textcolor{red}{(} g^{N}\circ f_{12}\textcolor{red}{)} (z_{0}))>0.
$
Here, $f_{12}(-1)$ is on $HL(\alpha,\textcolor{red}{(} g\circ f_{22}\textcolor{red}{)} (z_{0}))$,
because
\begin{align*}
\frac{f_{12}(-1)-\alpha}{\textcolor{red}{(} g\circ f_{22}\textcolor{red}{)} (z_{0})-\alpha}
= \frac{1-|\alpha|^{4}}{|\alpha|^{2}}>0
\end{align*}
holds by \eqref{f12} and \eqref{f22}. Hence we have
\begin{equation}\label{ConeLower}
\angle (\textcolor{red}{(} g\circ f_{22}\textcolor{red}{)} (z_{0}),\alpha,
\textcolor{red}{(} g^{N}\circ f_{12}\textcolor{red}{)} (z_{0})) >0.
\end{equation}
\eqref{ConeUpper} and \eqref{ConeLower} show that
$\textcolor{red}{(} g^{N}\circ f_{12}\textcolor{red}{)} (z_{0})$ belongs to
${\rm Cone}(f_{22}(\widetilde{A_{1}}))\setminus HL(\alpha,\textcolor{red}{(} g\circ f_{22}\textcolor{red}{)} (z_{0}))$.
Thus $k_{1}$ in Lemma \ref{lemma10} can be chosen as $N$.

Let us show that $l_{1}$ in Lemma \ref{lemma10}
can be chosen as $N+4$. Let us consider
${\rm Cone}(\textcolor{red}{(} \psi\circ f_{22}\textcolor{red}{)} (\widetilde{A_{1}}))$, see Figure \ref{cones1}.
The above discussion between
${\rm Cone}(f_{12}(\widetilde{A_{1}}))$
and ${\rm Cone}(f_{22}(\widetilde{A_{1}}))$
\textcolor{red}{works on} the discussion between
${\rm Cone}(f_{22}(\widetilde{A_{1}}))$
and ${\rm Cone}(\textcolor{red}{(} \psi\circ f_{22}\textcolor{red}{)} (\widetilde{A_{1}}))$.
Hence we see that $\textcolor{red}{(} g^{N}\circ f_{22}\textcolor{red}{)} (z_{0})$ belongs to
${\rm Cone}(\textcolor{red}{(} \psi\circ f_{22}\textcolor{red}{)} (\widetilde{A_{1}}))
\setminus HL(\alpha,\textcolor{red}{(} g\circ\psi\circ  f_{22}\textcolor{red}{)} (z_{0}))$.
Since
$\angle (\textcolor{red}{(} \psi\circ f_{22}\textcolor{red}{)} (z_{0}),\alpha, f_{22}(z_{0}))
= \angle (f_{22}(z_{0}),\alpha, f_{12}(z_{0}))=\theta$,
it follows that
$\angle (f_{12}(z_{0}), \alpha, \textcolor{red}{(} \psi\circ f_{22}\textcolor{red}{)} (z_{0}))
=2\pi - 2\theta=4\xi$. This means that
$g^{4}({\rm Cone}(\textcolor{red}{(} \psi\circ f_{22}\textcolor{red}{)} (\widetilde{A_{1}})))$
is identical to ${\rm Cone}(f_{12}(\widetilde{A_{1}}))$.
Thus $l_{1}$ in Lemma \ref{lemma10}
can be chosen as $N+4$.
\qed
\vspace{3mm}

Now, we study the range of $\xi$ for which both conditions (i) and (ii)
of Lemma \ref{hidarigawa} are satisfied.
\vspace{3mm}

\noindent{\bf The case $N\ge 4$ :}
As in the proof of Lemma \ref{hidarigawa},
$\textcolor{red}{(} g^{N}\circ f_{12}\textcolor{red}{)} (z_{0})$ belongs to
${\rm Cone}(f_{22}(\widetilde{A_{1}}))$.
Assume that the inequality
$|\textcolor{red}{(} g^{N}\circ f_{12}\textcolor{red}{)} (z_{0})-\alpha|
< |\textcolor{red}{(} g\circ f_{2212}\textcolor{red}{)} (z_{0})-\alpha|$ is satisfied.
Then $\textcolor{red}{(} g^{N}\circ f_{12}\textcolor{red}{)} (z_{0})$
belongs to the open ball
$B_{|\textcolor{red}{(} g\circ f_{2212}\textcolor{red}{)} (z_{0})-\alpha|}^{\circ}(\alpha)$.
These and
$\angle(\alpha,\textcolor{red}{(} g\circ f_{2212}\textcolor{red}{)} (z_{0}),f_{2212}(z_{0}))=\theta > \pi/2$
give
\begin{align*}
\textcolor{red}{(} g^{N}\circ f_{12}\textcolor{red}{)} (z_{0})
& \in {\rm Cone}(f_{22}(\widetilde{A_{1}})) \cap
B_{|\textcolor{red}{(} g\circ f_{2212}\textcolor{red}{)} (z_{0})-\alpha|}^{\circ}(\alpha) \\
& \subset
\mathcal{P}(\alpha,f_{2212}(z_{0}),\textcolor{red}{(} g\circ f_{2212}\textcolor{red}{)} (z_{0}))
\setminus [f_{2212}(z_{0}),\textcolor{red}{(} g\circ f_{2212}\textcolor{red}{)} (z_{0})] \\
& \subset V^{+}(\textcolor{red}{(} g\circ f_{2212}\textcolor{red}{)} (z_{0}),f_{2212}(z_{0})).
\end{align*}
Thus, if the inequality
$|\textcolor{red}{(} g^{N}\circ f_{12}\textcolor{red}{)} (z_{0})-\alpha|
< |\textcolor{red}{(} g\circ f_{2212}\textcolor{red}{)} (z_{0})-\alpha|$
is satisfied, then the condition (i) in Lemma \ref{hidarigawa} is satisfied.

We prove that the inequality
$|\textcolor{red}{(} g^{N}\circ f_{12}\textcolor{red}{)} (z_{0})-\alpha|
< |\textcolor{red}{(} g\circ f_{2212}\textcolor{red}{)} (z_{0})-\alpha|$
is valid for $N\ge 4$.
By \eqref{f12}, \eqref{f22}, and the expression of $g$,
it is verified that
\begin{equation*}
\frac{\textcolor{red}{(}g^{N}\circ f_{12}\textcolor{red}{)}(z_{0})-\alpha}
{\textcolor{red}{(}g\circ f_{2212}\textcolor{red}{)}(z_{0})-\alpha}
=-\frac{\alpha}{\overline{\alpha}}\times
\frac{\alpha^{N-1}}{1-|\alpha|^{2}-|\alpha|^{4}}.
\end{equation*}
The range $\pi/(N+2)\le \xi< \pi/(N+1)$ gives
$\sqrt{2}< 2\cos(\pi/(N+1))< 2\cos \xi \le 2\cos(\pi/(N+2)) <2$.
Then, by $|\alpha|=1/(2\cos \xi)$ and $N\ge 4$,
\begin{align*}
& \Big|\frac{\textcolor{red}{(}g^{N}\circ f_{12}\textcolor{red}{)}(z_{0})-\alpha}
{\textcolor{red}{(} g\circ f_{2212}\textcolor{red}{)} (z_{0})-\alpha}\Big|
=\frac{1}{1-(2\cos \xi)^{-2}-(2\cos \xi)^{-4}}
\times \frac{1}{(2\cos \xi)^{N-1}} \\
& < \frac{1}{1-(2\cos(\pi/(N+1)))^{-2}-(2\cos(\pi/(N+1)))^{-4}}
\times \frac{1}{(2\cos(\pi/(N+1)))^{3}} \\
& \le \frac{1}{(2\cos(\pi/5))^{3}-2\cos(\pi/5)-\frac{1}{2\cos(\pi/5)}}
= \frac{1}{2} < 1.
\end{align*}
Thus the condition (i) in Lemma \ref{hidarigawa} is satisfied for $N\ge 4$.

As in the proof of Lemma \ref{hidarigawa},
$\textcolor{red}{(}g^{N+4}\circ f_{22}\textcolor{red}{)}(z_{0})$ belongs to
${\rm Cone}(f_{12}(\widetilde{A_{1}}))$.
Assume that the inequality
$|\textcolor{red}{(}g^{N+4}\circ f_{22}\textcolor{red}{)}(z_{0})-\alpha|
< |\textcolor{red}{(}g\circ f_{1212}\textcolor{red}{)}(z_{0})-\alpha|$ is satisfied.
Then $\textcolor{red}{(}g^{N+4}\circ f_{22}\textcolor{red}{)}(z_{0})$
belongs to the open ball
$B_{|\textcolor{red}{(} g\circ f_{1212}\textcolor{red}{)} (z_{0})-\alpha|}^{\circ}(\alpha)$.
These and
$\angle(\alpha,\textcolor{red}{(} g\circ f_{1212}\textcolor{red}{)} (z_{0}),f_{1212}(z_{0}))=\theta > \pi/2$
give
\begin{align*}
\textcolor{red}{(} g^{N+4}\circ f_{22}\textcolor{red}{)} (z_{0})
& \in {\rm Cone}(f_{12}(\widetilde{A_{1}})) \cap
B_{|\textcolor{red}{(} g\circ f_{1212}\textcolor{red}{)} (z_{0})-\alpha|}^{\circ}(\alpha) \\
& \subset
\mathcal{P}(\alpha,f_{1212}(z_{0}),\textcolor{red}{(} g\circ f_{1212}\textcolor{red}{)} (z_{0}))
\setminus [f_{1212}(z_{0}),\textcolor{red}{(} g\circ f_{1212}\textcolor{red}{)} (z_{0})] \\
& \subset V^{+}(\textcolor{red}{(} g\circ f_{1212}\textcolor{red}{)} (z_{0}),f_{1212}(z_{0})).
\end{align*}
Thus, if the inequality
$|\textcolor{red}{(} g^{N+4}\circ f_{22}\textcolor{red}{)} (z_{0})-\alpha|
< |\textcolor{red}{(} g\circ f_{1212}\textcolor{red}{)} (z_{0})-\alpha|$
is satisfied, then the condition (ii) in Lemma \ref{hidarigawa} is satisfied.

The inequality
$|\textcolor{red}{(} g^{N+4}\circ f_{22}\textcolor{red}{)} (z_{0})-\alpha|
< |\textcolor{red}{(} g\circ f_{1212}\textcolor{red}{)} (z_{0})-\alpha|$
is valid for $N\ge 4$, and so is even if $N=3$.
In fact, for any $N\ge 3$, it is similarly verified as above that
\begin{equation*}
\frac{\textcolor{red}{(} g^{N+4}\circ f_{22}\textcolor{red}{)} (z_{0})-\alpha}
{\textcolor{red}{(} g\circ f_{1212}\textcolor{red}{)} (z_{0})-\alpha}
=-\frac{|\alpha|^{2}\alpha^{N+1}}{1-|\alpha|^{2}-|\alpha|^{4}}
\end{equation*}
and
\begin{align*}
& \Big|\frac{\textcolor{red}{(} g^{N+4}\circ f_{22}\textcolor{red}{)} (z_{0})-\alpha}
{\textcolor{red}{(} g\circ f_{1212}\textcolor{red}{)} (z_{0})-\alpha}\Big|
=\frac{1}{1-(2\cos \xi)^{-2}-(2\cos \xi)^{-4}}
\times \frac{1}{(2\cos \xi)^{N+3}} \\
& < \frac{1}{1-(2\cos(\pi/(N+1)))^{-2}-(2\cos(\pi/(N+1)))^{-4}}
\times \frac{1}{(2\cos(\pi/(N+1)))^{6}} \\
& \le \frac{1}{(2\cos(\pi/4))^{6}-(2\cos(\pi/4))^{4}-(2\cos(\pi/4))^{2}}
=\frac{1}{2} < 1.
\end{align*}
Thus the condition (ii) in Lemma \ref{hidarigawa} is satisfied for $N\ge 3$.

We have seen that both conditions (i) and (ii) with $N\ge 4$
in Lemma \ref{hidarigawa} are satisfied. Hence
$f_{1}(\mathcal{C})\cap f_{2}(\mathcal{C})=\emptyset$
is satisfied in this case.
\vspace{3mm}

\noindent{\bf The case $N= 3$ :}
We have already seen in the previous case
that the condition (ii) with $N=3$ in Lemma \ref{hidarigawa} is satisfied.

We determine an angle $\xi_{0}$ so that
the condition (i) with $N=3$ in Lemma \ref{hidarigawa}
is satisfied for any $\xi$ in $\pi/5\le \xi< \xi_{0}< \pi/4$.
The condition
$\textcolor{red}{(} g^{3}\circ f_{12}\textcolor{red}{)} (z_{0})
\in V^{+}(\textcolor{red}{(} g\circ f_{2212}\textcolor{red}{)} (z_{0}),f_{2212}(z_{0}))$
is satisfied if and only if
\begin{equation}\label{imaginal}
\Im\frac{\textcolor{red}{(} g^{3}\circ f_{12}\textcolor{red}{)} (z_{0})
-\textcolor{red}{(} g\circ f_{2212}\textcolor{red}{)} (z_{0})}
{f_{2212}(z_{0})-\textcolor{red}{(} g\circ f_{2212}\textcolor{red}{)} (z_{0})}>0
\end{equation}
is satisfied.
By \eqref{f12}, \eqref{f22}, and the expression of $g$,
it is verified that
\begin{equation*}
\frac{\textcolor{red}{(} g^{3}\circ f_{12}\textcolor{red}{)} (z_{0})
-\textcolor{red}{(} g\circ f_{2212}\textcolor{red}{)} (z_{0})}
{f_{2212}(z_{0})-\textcolor{red}{(} g\circ f_{2212}\textcolor{red}{)} (z_{0})}
=-\frac{\alpha}{\overline{\alpha}}
+\Big(\frac{\alpha}{\overline{\alpha}}\Big)^{3}
-\Big(\frac{\alpha}{\overline{\alpha}}\Big)^{3}\times
\frac{1-|\alpha|^{4}}{1-|\alpha|^{2}-|\alpha|^{4}}.
\end{equation*}
Substituting
$-\alpha/\overline{\alpha}=-\cos (2\xi)+i \sin (2\xi)$
and
$(\alpha/\overline{\alpha})^{3}=\cos (6\xi)-i \sin (6\xi)$
into the right-hand side of the above, we have
\begin{equation*}
\Im\frac{\textcolor{red}{(} g^{3}\circ f_{12}\textcolor{red}{)} (z_{0})
-\textcolor{red}{(} g\circ f_{2212}\textcolor{red}{)} (z_{0})}
{f_{2212}(z_{0})-\textcolor{red}{(} g\circ f_{2212}\textcolor{red}{)} (z_{0})}
=\sin(2\xi)+
\frac{|\alpha|^{2}}{1-|\alpha|^{2}-|\alpha|^{4}}\times \sin(6\xi).
\end{equation*}
Since
$\sin (2\xi)=2\sin\xi \cos\xi$,
$\sin (6\xi)=32\sin\xi \cos^{5}\xi -32\sin\xi \cos^{3}\xi +6\sin\xi \cos\xi$,
$|\alpha|=1/x$, and $x=2\cos\xi$, we have
\begin{equation*}
\Im\frac{\textcolor{red}{(} g^{3}\circ f_{12}\textcolor{red}{)} (z_{0})
-\textcolor{red}{(} g\circ f_{2212}\textcolor{red}{)} (z_{0})}
{f_{2212}(z_{0})-\textcolor{red}{(} g\circ f_{2212}\textcolor{red}{)} (z_{0})}
=(\sin\xi)\times \frac{x}{x^{4}-x^{2}-1}\times
(x^{6}-3x^{4}+2x^{2}-1).
\end{equation*}
The range $\pi/5\le \xi <\pi/4$ gives $\sqrt{2}<x\le (1+\sqrt{5})/2$.
For $\sqrt{2}<x\le (1+\sqrt{5})/2$,
$\sin\xi>0$ and $x/(x^{4}-x^{2}-1)>0$ hold.
The polynomial $x^{6}-3x^{4}+2x^{2}-1$ takes
a negative (resp.\ positive) value at $x=\sqrt{2}$
(resp.\ $x=(1+\sqrt{5})/2$), and
it is strictly monotonically increasing
for $\sqrt{2}<x\le (1+\sqrt{5})/2$.
Hence, there exists only one solution $x_{0}$
of $x^{6}-3x^{4}+2x^{2}-1=0$, whose approximate value
is $x_{0}\approx 1.5247$,
such that $\sqrt{2}<x_{0}< (1+\sqrt{5})/2$
and $x^{6}-3x^{4}+2x^{2}-1>0$ for $x_{0}<x\le (1+\sqrt{5})/2$.
Thus \eqref{imaginal} is satisfied for $x_{0}<x\le (1+\sqrt{5})/2$.
Put $\xi_{0}=\arccos (x_{0}/2)$, whose approximate value
is $\xi_{0}\approx 0.703858$. Then
$\textcolor{red}{(} g^{3}\circ f_{12}\textcolor{red}{)} (z_{0})
\in V^{+}(\textcolor{red}{(} g\circ f_{2212}\textcolor{red}{)} (z_{0}),f_{2212}(z_{0}))$
is satisfied for any $\xi$ in $\pi/5\le \xi< \xi_{0}< \pi/4$.
\vspace{3mm}

We have seen that both conditions (i) and (ii) with $N=3$ in Lemma \ref{hidarigawa}
are satisfied for any $\xi$ in $\pi/5\le \xi< \xi_{0}< \pi/4$.
Hence $f_{1}(\mathcal{C})\cap f_{2}(\mathcal{C})=\emptyset$
is satisfied for any $\xi$ in this range.
By combining the result of the case $N=3$ with that of the case $N\ge 4$,
the proof of Proposition \ref{prop4} is completed.

\section{Completion of proof of Theorem \ref{MainTheorem}}\label{proofTh1}

We will prove that the set of limit points of $\mathcal{C}$
outside $\mathcal{C}$ is $\{0,1\}$, which comes from the sets
$\bigcup_{n=1}^{\infty}A_{n}$ and
$\bigcup_{n=1}^{\infty}f_{2}(A_{n})$ in $\mathcal{C}$.
Since $\bigcup_{n=1}^{\infty}f_{2}(A_{n})$ is the image of
$\bigcup_{n=1}^{\infty}A_{n}$ by $f_{2}$,
it \textcolor{red}{suffices} to prove that
the set of limit points of $\bigcup_{n=1}^{\infty}A_{n}$ outside
$\bigcup_{n=1}^{\infty}A_{n}$ is $\{0\}$.

Since $A_{n}\subset\widetilde{A_{n}}=
\mathcal{P}(0,f_{(1)^{n-1}}(z_{0}),f_{(1)^{n}}(z_{0}))$
and the sequence $(|f_{(1)^{n}}(z_{0})|)_{n=1}^{\infty}$
is decreasing as $n\to\infty$,
$A_{n}\subset B_{|f_{(1)^{m}}(z_{0})|}(0)$ holds
for all $n$ with $n\ge m+1$, where $m$ is a non-negative integer.
For any $\delta>0$, choose $m$ as
$m>\max\{\frac{\log\delta^{-1}}{\log \sqrt{2}},0\}$.
Then, by \eqref{z0} and $1/2 < |\alpha| < 1/\sqrt{2}$, we have
$
|f_{(1)^{m}}(z_{0}))|=|\alpha|^{m+1}/(1-|\alpha|^{4})
<(1/\sqrt{2})^{m}<\delta
$.
Hence $A_{n}\subset B_{\delta}(0)$ holds for all $n$ with
$n>\max\{\frac{\log\delta^{-1}}{\log \sqrt{2}},0\}+1$.

Next, let us show the statement: if $z$ is a limit point of
$\bigcup_{n=1}^{\infty}A_{n}$ and $z$ is not equal to $0$, then
$z$ belongs to $\bigcup_{n=1}^{\infty}A_{n}$.
By the assumption,
$B_{\varepsilon}(z)\cap \big(\bigcup_{n=1}^{\infty}A_{n}\big)
\neq \emptyset$ is satisfied for any $\varepsilon$ with $0<\varepsilon <|z|/3$.
Obviously,
$B_{\varepsilon}(z)\cap B_{|z|/3}(0)= \emptyset$.
By choosing the above $\delta$ as $\delta=|z|/3$,
$A_{n}\subset B_{|z|/3}(0)$ holds for all $n$ with
$n>\max\{\frac{\log (3|z|^{-1})}{\log \sqrt{2}},0\}+1$.
Hence,
$B_{\varepsilon}(z)\cap \big(\bigcup_{n=1}^{\infty}A_{n}\big)
\neq \emptyset$ implies that for any $\varepsilon$ with $0<\varepsilon <|z|/3$
there exists a point $w_{\varepsilon}$ such that
\begin{equation}\label{tenretsu}
B_{\varepsilon}(z)\cap \mathcal{A}_{z}\ni w_{\varepsilon},\quad
\mbox{where
$\mathcal{A}_{z}
=\bigcup_{1\le n\le \max\{\frac{\log (3|z|^{-1})}{\log \sqrt{2}},0\}+1}A_{n}$.}
\end{equation}
Here $\mathcal{A}_{z}$ is compact,
because it is a finite union of the compact sets $A_{n}$,
and, moreover, it does not depend on $\varepsilon$. \eqref{tenretsu} means that
the sequence $(w_{\varepsilon})$ belongs to
the compact set $\mathcal{A}_{z}$ and it converges to $z$ as $\varepsilon\to 0$.
Therefore $z$ belongs to $\mathcal{A}_{z}$,
and, consequently, $z$ belongs to $\bigcup_{n=1}^{\infty}A_{n}$.

Finally, let us show the statement: $0$ is a limit point of
$\bigcup_{n=1}^{\infty}A_{n}$. The property that
$A_{n}\subset B_{\delta}(0)$ holds for any $\delta>0$ and all $n$
with $n>\max\{\frac{\log\delta^{-1}}{\log \sqrt{2}},0\}+1$ implies that
$B_{\delta}(0)\cap \big(\bigcup_{n=1}^{\infty}A_{n}\big)
\neq \emptyset$ is satisfied for any $\delta>0$, and, consequently,
$0$ is a limit point of $\bigcup_{n=1}^{\infty}A_{n}$.

We have seen that the set of limit points of $\mathcal{C}$ outside $\mathcal{C}$
is $\{0,1\}$, i.e., $\overline{\mathcal{C}}=\mathcal{C}\cup \{0,1\}$.
This and Proposition \ref{prop4} give that, for $0<\xi<\xi_{0}$,
\begin{equation}\label{property3}
	f_{1}(\overline{\mathcal{C}})\cap f_{2}(\overline{\mathcal{C}})
	=(f_{1}(\mathcal{C})\cup \{0,\alpha\})
	\cap (f_{2}(\mathcal{C}) \cup \{1,\alpha\})=\{\alpha\}.
\end{equation}

We now choose the open set $U_{\xi}$ to be the interior of $\mathcal{C}$.
Then, by Propositions \ref{prop2} and \ref{prop4}, and by \eqref{property3},
for $0<\xi<\xi_{0}$,
we have $U_{\xi}\supset f_{1} (U_{\xi}) \cup f_{2}(U_{\xi})$,
$f_{1}(U_{\xi})\cap f_{2}(U_{\xi})=\emptyset$,
and $f_{1}(\overline{U_{\xi}})\cap f_{2}(\overline{U_{\xi}})=\{\alpha\}$.
Hence, by Proposition \ref{Arc}, we conclude that
the limit Dragon curve ${\mathcal D}(\xi)$ is a simple arc for $0<\xi<\xi_{0}$.
This completes the proof of Theorem \ref{MainTheorem}.

\section{Proof of Theorem \ref{Main2}}\label{proofTh2}

To prove Theorem \ref{Main2}, we introduce a finite polygon $C_{k}$,
which covers the renormalized \textcolor{red}{Dragon curve of order $k$, $D_{k}$,}
in Definition \ref{def1} except for its beginning and end segments.
$C_{k}$ is regarded as a truncated set of $\mathcal{C}$ of Definition \ref{def4}.

\begin{figure}
\begin{center}
\caption{The renormalized \textcolor{red}{Dragon curve of order $5$} and the set $C_{5}$}
\label{truncate}
\includegraphics[clip, width=7.5cm]{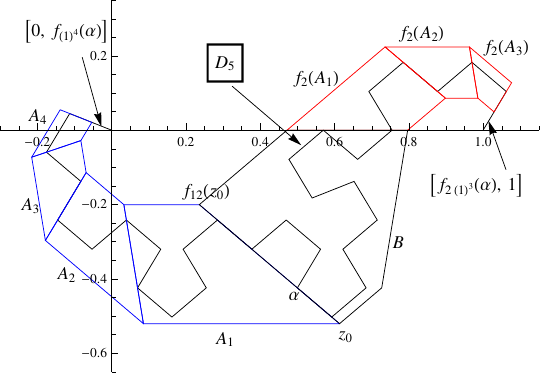}
\end{center}
\end{figure}

\begin{defn}\label{def6}
Under the same \textcolor{red}{notation} as in Definition \ref{def4},
the set $C_{k}$, $k\ge 2$, is defined by
$$
C_{k}=\Big(\bigcup_{n=1}^{k-1}A_{n}\Big)
\cup B \cup \Big(\bigcup_{n=1}^{k-2}f_{2}(A_{n})\Big),
$$
where $\bigcup_{n=1}^{0}f_{2}(A_{n})$ is defined to be $\emptyset$,
see Figure \ref{truncate}.
\end{defn}

\begin{figure}
\begin{center}
\caption{The sets $f_{1}(C_{5})$ and $f_{2}(C_{5})$}
\label{nointersection}
\includegraphics[clip, width=7.5cm]{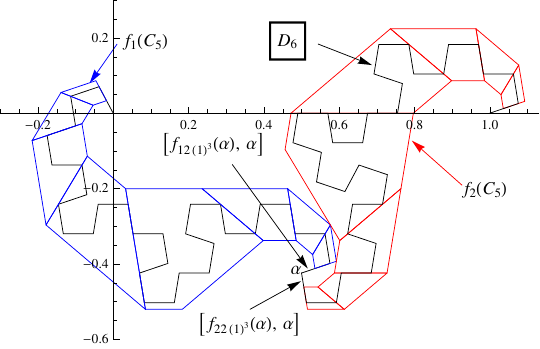}
\end{center}
\end{figure}

\begin{lemma}\label{lemmaT1}
Let $\xi$ be a fixed angle in $0<\xi<\pi/4$, $k\ge 2$,
and $C_{k}$ be the set of Definition \ref{def6}.
Then {\rm (i)} $f_{1}(C_{k})\subset C_{k}\cup A_{k}$ and
{\rm (ii)} $f_{2}(C_{k})\subset C_{k}\cup f_{2}(A_{k-1})$
are satisfied, see Figure \ref{nointersection}.
\end{lemma}

\noindent{\it Proof.}
(i) \textcolor{red}{In} the expression
$
f_{1}(C_{k})=\Big(\bigcup_{n=2}^{k-1}A_{n}\Big)\cup A_{k}
\cup f_{1}(B) \cup \Big(\bigcup_{n=1}^{k-2}f_{12}(A_{n})\Big)
$,
as in the proof of Proposition \ref{prop2}, we see that
$f_{1}(B)\subset A_{1}$ and
$\bigcup_{n=1}^{k-2}f_{12}(A_{n})\subset \widetilde{B}\subset A_{1}\cup B$.
Hence we obtain the desired inclusion.

(ii) \textcolor{red}{In} the expression
$
f_{2}(C_{k})=\Big(\bigcup_{n=1}^{k-2}f_{2}(A_{n})\Big)\cup f_{2}(A_{k-1})
\cup f_{2}(B) \cup \Big(\bigcup_{n=1}^{k-2}f_{22}(A_{n})\Big)
$,
as in the proof of Proposition \ref{prop2}, we see that
$f_{2}(B)\subset \widetilde{B}\subset A_{1}\cup B$ and
$\bigcup_{n=1}^{k-2}f_{22}(A_{n})\subset \widetilde{B}\subset A_{1}\cup B$.
Hence we obtain the desired inclusion.
\qed

\begin{lemma}\label{lemmaT2}
Let $\xi$ be a fixed angle in $0<\xi<\pi/4$, $k\ge 2$,
$D_{k}$ be the renormalized \textcolor{red}{Dragon curve of order $k$},
and $C_{k}$ be the sets of Definition \ref{def6}. Then
$
D_{k}\subset C_{k}\cup [0,f_{(1)^{k-1}}(\alpha)]\cup [f_{2(1)^{k-2}}(\alpha),1]
$
holds.
\end{lemma}

\noindent{\it Proof.}
The statement is proved by induction on $k\ge 2$.
$D_{2}$ is the broken line such that
$
D_{2}=[0,f_{1}(\alpha)]\cup [f_{1}(\alpha),\alpha]
\cup [\alpha,f_{2}(\alpha)]\cup [f_{2}(\alpha),1].
$
It is verified that $\alpha$ and $f_{1}(\alpha)$ are on the boundary of $A_{1}$,
and $\alpha$ and $f_{2}(\alpha)$ are on the boundary of $B$.
Hence,
$
[f_{1}(\alpha),\alpha]\cup [\alpha,f_{2}(\alpha)]
\subset A_{1}\cup B=C_{2},
$
and the statement for $k=2$ is valid.

Assume that the statement for $k$ is satisfied. Then,
$
f_{1}(D_{k})\subset f_{1}(C_{k})\cup [0,f_{(1)^{k}}(\alpha)]
\cup [f_{12(1)^{k-2}}(\alpha),\alpha]
$
and
$
f_{2}(D_{k})\subset f_{2}(C_{k})\cup [1,f_{2(1)^{k-1}}(\alpha)]
\cup [f_{22(1)^{k-2}}(\alpha),\alpha]
$.
These inclusions, the definition of $D_{k}$, and Lemma \ref{lemmaT1} give
\begin{align*}
D_{k+1}& =f_{1}(D_{k}) \cup f_{2}(D_{k}) \\
& \subset C_{k}\cup A_{k}\cup f_{2}(A_{k-1})
\cup [f_{12(1)^{k-2}}(\alpha),\alpha]
\cup [f_{22(1)^{k-2}}(\alpha),\alpha] \\
& \qquad \cup [0,f_{(1)^{k}}(\alpha)]\cup [1,f_{2(1)^{k-1}}(\alpha)].
\end{align*}
By the definition of $C_{k}$, it holds that
$C_{k}\cup A_{k}\cup f_{2}(A_{k-1})=C_{k+1}$.
For the set $\widetilde{B}$ in the proof of Proposition \ref{prop2},
we have $\widetilde{B}\subset A_{1}\cup B=C_{2}\subset C_{k+1}$.
Thus it \textcolor{red}{suffices} to prove that both
$[f_{12(1)^{k-2}}(\alpha),\alpha]$ and $[f_{22(1)^{k-2}}(\alpha),\alpha]$
are included in $\widetilde{B}$.
It is verified that $[f_{(1)^{k-2}}(\alpha),0]\subset \widetilde{A_{k-1}}$.
This and
$\bigcup_{n=1}^{\infty}\widetilde{A_{n}}\subset S'$ (see Lemma \ref{lemma3})
give $[f_{(1)^{k-2}}(\alpha),0]\subset S'$.
Applying $f_{12}$ (resp.\ $f_{22}$) to $[f_{(1)^{k-2}}(\alpha),0]\subset S'$,
then using $f_{12}(S')\subset \widetilde{B}$
(resp.\ $f_{22}(S')\subset \widetilde{B}$),
which has been shown in the proof of Proposition \ref{prop2},
we have $[f_{12(1)^{k-2}}(\alpha),\alpha]\subset \widetilde{B}$
(resp.\ $[f_{22(1)^{k-2}}(\alpha),\alpha]\subset \widetilde{B}$).
\qed

\begin{lemma}\label{lemmaT3}
Let $\xi$ be a fixed angle in $0<\xi<\pi/4$, $k\ge 2$,
and $D_{k}$ be the renormalized \textcolor{red}{Dragon curve of order $k$}.
Under the same \textcolor{red}{notation} as in Lemma \ref{lemma5},
we have the \textcolor{red}{following};

{\rm (i)} $D_{k}\cap \psi(S)
\subset \big(\bigcup_{m=-1}^{k-2}f_{2}(A_{m})\big)\cup [f_{2(1)^{k-2}}(\alpha),1]$,

{\rm (ii)} $D_{k}\cap \tau(S)
\subset \big(\bigcup_{m=-1}^{k-2}f_{2}(A_{m})\big)\cup [f_{2(1)^{k-2}}(\alpha),1]$.
\end{lemma}

\noindent{\it Proof.}
(i) By Lemma \ref{lemmaT2} and $[0,f_{(1)^{k-1}}(\alpha)]\subset \widetilde{A_{k}}$,
$$
D_{k}\cap \psi(S)
\subset (C_{k}\cap \psi(S))\cup (\widetilde{A_{k}}\cap \psi(S))
\cup ([f_{2(1)^{k-2}}(\alpha),1]\cap \psi(S)).
$$
\textcolor{red}{In} the expression
$
C_{k}=\big(\bigcup_{n=2}^{k-1}A_{n}\big)
\cup (A_{1}\cup B) \cup \big(\bigcup_{n=1}^{k-2}f_{2}(A_{n})\big),
$
each of the three sets on the right-hand side is treated by using
(i), (ii), and (iii) in the proof of Lemma \ref{lemma7}.
Then we have
$
C_{k}\cap \psi(S)\subset \bigcup_{m=-1}^{k-2}f_{2}(A_{m}).
$
The discussion for (iii) in the proof of Lemma \ref{lemma7} can
be modified to get
$\big(\bigcup_{n=2}^{\infty}\widetilde{A_{n}}\big) \cap \psi(S)=\emptyset$,
and especially, $\widetilde{A_{k}}\cap \psi(S)=\emptyset$.
The discussion for (i) in the proof of Lemma \ref{lemma7} can
be modified to get
$\bigcup_{n=1}^{\infty}f_{2}(\widetilde{A_{n}}) \subset \psi(S)$,
and especially, $f_{2}(\widetilde{A_{k-1}})\subset \psi(S)$.
Since $[f_{2(1)^{k-2}}(\alpha),1]\subset f_{2}(\widetilde{A_{k-1}})$,
we have $[f_{2(1)^{k-2}}(\alpha),1]\cap \psi(S)=[f_{2(1)^{k-2}}(\alpha),1]$.
Thus the desired inclusion is obtained.

(ii) By Lemmas \ref{lemma3} and \ref{lemma5}, we see that
$C_{k}\subset \mathcal{C}\subset T'$,
$[0,f_{(1)^{k-1}}(\alpha)]\subset \widetilde{A_{k}}\subset S\subset T'$,
and $[f_{2(1)^{k-2}}(\alpha),1]\subset \psi(\widetilde{A_{k}})
\subset \psi(S)\subset T'$. Combining these with Lemma \ref{lemmaT2}, we have
$D_{k}\subset T'$. By $D_{k}\subset T'$ and
the same argument as in the proof of Lemma \ref{lemma8},
$D_{k}\cap \tau(S) \subset D_{k}\cap \psi(S)$ is obtained.
Thus (ii) follows from (i). \qed

\begin{lemma}\label{lemmaT4}
Let $\xi$ be a fixed angle in $0<\xi<\pi/4$, $k\ge 2$,
and $D_{k}$ be the renormalized \textcolor{red}{Dragon curve of order $k$}. If
\begin{equation}\label{majiwaranaiT}
\Big(\Big(\bigcup_{m=-1}^{k-2}f_{12}(A_{m})\Big)
\cup [f_{12(1)^{k-2}}(\alpha),\alpha]\Big)
\cap
\Big(\Big(\bigcup_{m=-1}^{k-2}f_{22}(A_{m})\Big)
\cup [f_{22(1)^{k-2}}(\alpha),\alpha]\Big)
=\{\alpha\}
\end{equation}
is satisfied, then
$f_{1}(D_{k})\cap f_{2}(D_{k})=\{\alpha\}$ is satisfied.
\end{lemma}

\noindent{\it Proof.} This lemma is obtained by
Lemma \ref{lemmaT3}, Lemma \ref{lemma6}, and
the same argument as in the proof of Proposition \ref{prop3}. \qed

\vspace{3mm}

As a result of Section \ref{proofII},
the condition \eqref{majiwaranai} is satisfied for $0< \xi< \xi_{0}$.
Hence $\big(\bigcup_{m=-1}^{k-2}f_{12}(A_{m})\big)
\cap \big(\bigcup_{m=-1}^{k-2}f_{22}(A_{m})\big)=\emptyset$ is satisfied
for $0< \xi< \xi_{0}$.
Since the segments $[f_{12(1)^{k-2}}(\alpha),\alpha]$ and
$[f_{22(1)^{k-2}}(\alpha),\alpha]$ make the angle $\theta$,
$[f_{12(1)^{k-2}}(\alpha),\alpha]\cap [f_{22(1)^{k-2}}(\alpha),\alpha]=\{\alpha\}$
holds. Hence the condition \eqref{majiwaranaiT} is satisfied
for $0< \xi< \xi_{0}$, if both conditions
\begin{align}
& [f_{12(1)^{k-2}}(\alpha),\alpha]\cap \Big(\bigcup_{m=-1}^{k-2}f_{22}(A_{m})\Big)
=\emptyset, \quad 0< \xi< \xi_{0}, \label{f22spiral} \\
& \Big(\bigcup_{m=-1}^{k-2}f_{12}(A_{m})\Big)
\cap [f_{22(1)^{k-2}}(\alpha),\alpha]=\emptyset, \quad 0< \xi< \xi_{0},
\label{f12spiral}
\end{align}
are satisfied.
By Lemma \ref{lemma9} and
$[f_{12(1)^{k-2}}(\alpha),\alpha]\subset f_{12}(\widetilde{A_{k-1}})
=\textcolor{red}{(} g^{k-2}\circ f_{12}\textcolor{red}{)} (\widetilde{A_{1}})$
(resp.\
$[f_{22(1)^{k-2}}(\alpha),\alpha]\subset f_{22}(\widetilde{A_{k-1}})
=\textcolor{red}{(} g^{k-2}\circ f_{22}\textcolor{red}{)} (\widetilde{A_{1}})$),
the set on the left-hand side of \eqref{f22spiral} (resp.\ \eqref{f12spiral})
is contained in
$
\textcolor{red}{(} g^{k-2}\circ f_{12}\textcolor{red}{)} (\widetilde{A_{1}})\cap
\big(\bigcup_{m=-2}^{k-3}\textcolor{red}{(} g^{m}\circ f_{22}\textcolor{red}{)} (A_{1})\big)
$
(resp.\
$
\big(\bigcup_{m=-2}^{k-3}\textcolor{red}{(} g^{m}\circ f_{12}\textcolor{red}{)} (A_{1})\big)
\cap \textcolor{red}{(} g^{k-2}\circ f_{22}\textcolor{red}{)} (\widetilde{A_{1}})$).
Then, by the same arguments as in the beginning of Section \ref{proofII},
we see that if both conditions
\begin{equation}\label{emptyTA}
\Big(\bigcup_{n=1}^{k}\textcolor{red}{(} g^{n}\circ f_{12}\textcolor{red}{)} (\widetilde{A_{1}})\Big)
\cap {\rm Cone}(f_{22}(\widetilde{A_{1}}))
\cap f_{22}(A_{1})=\emptyset
\end{equation}
and
\begin{equation}\label{emptyTB}
f_{12}(A_{1})\cap
\Big(\bigcup_{n=1}^{k}\textcolor{red}{(} g^{n}\circ f_{22}\textcolor{red}{)} (\widetilde{A_{1}})\Big)
\cap {\rm Cone}(f_{12}(\widetilde{A_{1}})) =\emptyset
\end{equation}
are satisfied, then both conditions \eqref{f22spiral} and \eqref{f12spiral}
are satisfied.
Indeed, \eqref{emptyTA} and \eqref{emptyTB} are valid,
because \eqref{emptyA} and \eqref{emptyB} are valid for $0< \xi< \xi_{0}$
as a result of Section \ref{proofII}.
Thus, by Lemma \ref{lemmaT4}, we obtain
\begin{equation}\label{intersection}
f_{1}(D_{k})\cap f_{2}(D_{k})=\{\alpha\}, \quad 0< \xi< \xi_{0}.
\end{equation}

Let us prove Theorem \ref{Main2} by induction on $k\ge 1$.
$D_{1}$ is the broken line such that $D_{1}=[0,\alpha]\cup [\alpha,1]$.
Since the segments $[0,\alpha]$ and $[\alpha,1]$ make the angle $\theta$,
$D_{1}$ has no self-intersection.
If $D_{k}$ has no self-intersection,
then each of $f_{1}(D_{k})$ and $f_{2}(D_{k})$ has no self-intersection.
By these properties and the expression $D_{k+1}=f_{1}(D_{k}) \cup f_{2}(D_{k})$,
if $D_{k+1}$ has a self-intersection,
then it has to belong to $f_{1}(D_{k})\cap f_{2}(D_{k})$,
and moreover, it is $\alpha$ by \eqref{intersection}.
In $D_{k+1}$, both segments $[f_{12(1)^{k-2}}(\alpha),\alpha]$
and $[f_{22(1)^{k-2}}(\alpha),\alpha]$ contain $\alpha$.
Assume that there exists another segment $L\subset D_{k+1}$
such that $L$ contains $\alpha$. Without loss of generality,
we may assume that $L$ belongs to $f_{2}(D_{k})$. Then,
$f_{2}^{-1}(\alpha)\in [f_{2(1)^{k-2}}(\alpha),1]\cap f_{2}^{-1}(L)\subset D_{k}$,
hence $D_{k}$ has the self-intersection $f_{2}^{-1}(\alpha)$,
which is the point $1$.
This contradicts that $D_{k}$ has no self-intersection.
\textcolor{red}{Hence no segment except for $[f_{12(1)^{k-2}}(\alpha),\alpha]$
and $[f_{22(1)^{k-2}}(\alpha),\alpha]$ can contain $\alpha$.}
Thus $D_{k+1}$ has no self-intersection.

This completes the proof of Theorem \ref{Main2}.

\section{Future Problems}\label{Fin}

We enlist several open problems.

\begin{enumerate}
\item[(i)] Can we find the supremum of $\xi$, such that
$\mathcal{D}(\xi)$ is non\textcolor{red}{-}self-intersective?
The shape $\mathcal{C}$ can be modified to have a better constant,
but we do not know how to compute the supremum.
\item[(ii)]
We have seen that, for $0< \xi< \xi_{0}$, the IFS$\{f_1,f_2\}$
satisfies the open set condition.
Can we characterize when $\mathcal{D}(\xi)$
satisfies the open set condition?
Since non\textcolor{red}{-}self-intersectivity is stronger than open set condition,
we may have a strictly larger supremum than (i).
\item[(iii)] What can we say when $\xi>\pi/3$?
In this case, the system is expanding,
and it may be interesting to study the basic property of
point set $X=\{f_{\omega}(0) :
\omega=\omega_{1}\omega_{2}\omega_{3}\cdots, \omega_{j}\in\{1,2\} \}$.
Is $X$ uniformly discrete and/or relatively dense in ${\bf C}$?
\end{enumerate}

\end{document}